\documentclass[12pt]{article}
\makeatletter
\@addtoreset{equation}{section}
\makeatother

\def\sqr#1#2{{\vcenter{\vbox{\hrule height .#2pt
                             \hbox{\vrule width .#2pt height#1pt \kern#1pt
                                   \vrule width .#2pt}
                             \hrule height .#2pt}}}}
\def\square{\mathchoice\sqr54\sqr54\sqr{4.1}3\sqr{3.5}3}

\newtheorem{Th}{Theorem}[section]
\newtheorem{Pro}{Proposition}[section]
\newtheorem{Le}{Lemma}[section]
\newtheorem{Rem}{Remark}[section]
\newcommand{\be}[1]{\begin{equation}\label{#1}}
\newcommand{\ba}[1]{\begin{eqnarray}\label{#1}}
\newcommand{\ee}{\end{equation}}
\newcommand{\bl}[1]{\begin{Le}\label{#1}}
\newcommand{\bt}[1]{\begin{Th}\label{#1}}
\newcommand{\bp}[1]{\begin{Pro}\label{#1}}
\newcommand{\br}[1]{\begin{Rem}\label{#1}}
\def\z{\noindent}

\def\RR{\hbox{ I\hskip -2pt R}}

\def\NN{\hbox{\it I\hskip -2pt N}}

\def\z{\noindent}

\title{Optimal Uniform Estimates and Rigorous Asymptotics Beyond all Orders for
a Class of Ordinary Differential Equations
\thanks{Work supported, in part, by AFOSR Grant AF-0115}}
\author{Ovidiu Costin {\lowercase {and}} Martin D. Kruskal
\thanks {Mathematics Department, Hill Center
Rutgers University, New Brunswick, NJ 08903;
e-mail: costin@maxwell.rutgers.edu, kruskal@math.rutgers.edu}}
\begin{document}
\include{psfig}
\maketitle
\begin{abstract}
For first order differential equations of the form $y'=\sum_{p=0}^P
F_p(x)y^p$ and second order homogeneous linear differential equations
$y''+a(x)y'+b(x)y=0$ with locally integrable coefficients having
asymptotic (possibly divergent) power series when $|x|\rightarrow\infty$
on a ray $\arg(x)=$const,
under some further assumptions, it is shown that, on the given ray,
there is a one-to-one correspondence between true solutions
and (complete) formal solutions.  The correspondence is based on
asymptotic inequalities which are required to be uniform in $x$ and
optimal with respect to certain weights.
\end{abstract}

\begin{section}{Introduction and Main Results}

The main purpose of the present paper is to give, in terms of uniform
asymptotic estimates, a precise meaning to {\em complete}
asymptotic expansions (e.g., as power series followed by exponentially
small terms) of solutions of a class of differential equations in a
neighborhood of an irregular singular point (chosen to be infinity).
The study of the exponentially small terms in asymptotic
expansions has known a rapid development in the last years, 
especially after the pioneering works of Ecalle \cite{Ecalle-book}
and Berry \cite {Berry}.

For first order polynomially nonlinear (or linear) and second order
homogeneous linear ordinary differential equations we give, using
asymptotic inequalities, a one-to-one correspondence between formal
solutions and true solutions.  The representation of a given function
$y$ turns out to be (under some assumptions) independent of the
differential equation(s) of which $y$ is a solution.

The key ingredient is the concept of  uniform optimal asymptotic
inequalities which we illustrate in the following.

Consider a formal power series $\tilde
f=\sum_{k=0}^{\infty}f_k x^{-k}$, where $f_k$ are complex numbers
and $x$ is thought of being a large variable. We say that a
function $f$ is uniformly asymptotic to the series $\tilde f$ 
with respect to the weight $w$, along a given ray in the complex
plane, say ${\cal R}=x>x_0>0$ iff 

\be{op1}
\left|f(x)-\sum_{k=0}^{n-1} f_k x^{-k}\right|<w(n) x^{-n}
\ \ \forall x\in {\cal R}\ {\mbox{ and }} \forall n\in \NN
\ee

\z Any function $f$ that is asymptotic to the series $\tilde f$
is uniformly asymptotic with respect to {\em some} $w$. The minimal
such $w$ is obviously

\be{deftrans}
w(f,\tilde f;n)\equiv w_f(n):=\sup_{x>x_0}x^n\left|f(x)-\sum_{k=0}^{n-1}f_k
x^{-k}\right|
\ee

\z (we might think of $w_f$ as a transform of $f$ with respect 
to $\tilde f$). With no restrictions on $w$ we thus obtain the
Poincar\'e asymptoticity. With more conditions on $w$ we obtain
sharper asymptoticity classes: e.g. if $f_k\sim \alpha^n (n!)^{\beta}$
and we consider the functions $f$ for which $w_f(k)\le \gamma^n
(n!)^{\beta},\ \gamma>\alpha$ we obtain the familiar Gevrey or
Gevrey-Roumieux classes \cite{Ramis}.

It is natural to take $w_f$ as a measure of the separation between a
function $f$ and a formal series, $\tilde f$. We will say that $f$ is
closer to $\tilde f$ than $g$ iff

\be{deflessmax}
w_f(n)\le w_g(n)\  \mbox{ for all large $n$}
\ee
or, on occasions if a weaker condition holds:

\be{deflessdom}
\limsup_{n\rightarrow\infty}
\frac{w_f(n)}{ w_g(n)}\le 1
\ee

Given a formal series $\tilde f$, a ray ${\cal R}$
 and a class of functions ${\cal F}$,
it is also natural to ask what are the sharpest bounds compatible with it,
i.e. what is the ``greatest lower bound'' of the weights $w$ such that
$w=w_f$ for {\em some} $f$.

We say that a function $f$ is optimally  asymptotic  to $\tilde f$
with respect to 
${\cal F}$ along ${\cal R}$, written

\be{defappr}
f\simeq \sum_{k=0}^{\infty} f_k x^{-k} 
\mbox{ as } x\rightarrow\infty\ \ \ \mbox{on } {\cal R}
\ee

\z and correspondingly that
$w_f$ is an optimal weight for $\tilde f$
\z iff (\ref{deflessmax}) (or sometimes  (\ref{deflessdom})) holds 
for {\em all} $g\in{\cal F}$.

As a first example, note that if $\tilde f_0$ is a {\em convergent} series
and $x_0>\rho^{-1}>R^{-1}$ where $R$ is the radius of convergence,
and ${\cal F}$ is the set of all functions defined for $x>x_0$
then

$$f\simeq \tilde f_0 \Longleftrightarrow f=f_0=\sum_{k=0}^{\infty}f_{0;k}x^{-k}$$

\z Indeed,
$w_{f_0}(n)\le C\rho^n \left(1-\rho
x_0^{-1}\right)^{-1}\rightarrow 0$ as $n\rightarrow\infty$
and therefore uniform asymptoticity of a function $f$ with respect
to $w_0$ implies the convergence of the series $\tilde f_0$ to $f$.

Next, take 
$\tilde f_1=\sum_{k=0}^{\infty}k! x^{-k-1}$. We will see that
an optimal weight for this series,
in the sense (\ref{deflessdom}), along the ray $x>x_0>0$
behaves like $w(n)=(a_*\sqrt{n}+O(1))n!$ for large $n$, where $a_*$
is defined in (\ref{defA12}). An optimally asymptotic
function is $f(x)=e^{-x}Ei(x):=e^{-x} P\int_{-\infty}^x t^{-1}e^tdt$.

It turns out that if $f_1,\ f_2$ are 
two functions optimally asymptotic to $\tilde f_1$ then 
$f_1-f_2=o(e^{-x})$ for large $x$. 
It is then meaningful to extend the notion of optimal asymptoticity to
more general structures e.g., power series followed by exponentially
small terms: in our example we can {\em define}

\be{introexpsmall}
f(x)\simeq\sum_{k=0}^{\infty}\frac{k!}{ x^{k+1}}+Ce^{-x}
\Longleftrightarrow
f(x)-Ce^{-x}\simeq\sum_{k=0}^{\infty}\frac{k!}{ x^{k+1}}
\ee

\z since the second relation can hold for at most one value of $C$.
On the other hand $e^{-x}Ei(x)+Ce^{-x}$ is the general solution of the
equation $f'+f=x^{-1}$ whereas $\tilde f+Ce^{-x}$ is the general
formal solution of that equation. It follows that the relation
$f(x)\simeq\tilde f+C e^{-x}$ establishes a one-to -one correspondence
between true and formal solutions for this equation.

The main goal of the present paper is to show that by letting ${\cal
F}$ to contain solutions of a class of first and second order
differential equations, along some ray, optimal asymptoticity (in
either the stronger form (\ref{deflessmax}) or the weaker one
(\ref{deflessdom})) always gives a one-to one correspondence between
true and formal solutions.

As a manifestation
of the Stokes phenomenon, the complete expansion of a given solution
of a differential equation
will depend on the ray along which the asymptotic estimates are
considered. In our example, we have

\be{eqEi}
e^{-x}{\cal E}i(x)\simeq \sum_{j=0}^{\infty}\frac{j!}{x^{j+1}}
+\left\{\begin{array}{ll}
0&\ \ \mbox{if $\arg(x)=0$}\cr
\pi i\, e^{-x}&\ \ \mbox{if $\arg(x)\in (0,\pi)$}\cr
-\pi i\,e^{-x} &\ \ \mbox{if $\arg(x)\in (-\pi,0)$}
\end{array}\right.
\ee

For a complete study of the Stokes one would need to consider
besides rays, parabolas 
(of the form $\arg(x)=\lambda\,|x|^{-1/2}$
for (\ref{eqEi})); then the constant beyond all orders
changes smoothly in a narrow region near the Stokes line
(it is a smooth function of $\lambda $ in the example 
above), see \cite{Berry}. 

However in this paper we are not making analyticity assumptions on the
coefficients of the differential equations  and we will be consequently
mainly concerned with the behavior of solutions along a fixed ray.

A function is characterized by its generalized expansion
with a precision comparable to that of truncation to the least
term of the complete asymptotic expansion. The relation with this technique
is explored in section 2.

In this respect we mention the  paper \cite{McLeod} 
in which it
is shown that truncation to the least term can be used to measure
the terms beyond all orders for second order linear homogeneous
differential equations with coefficients analytic at infinity.
The
corresponding results in our paper show that the same is true for
equations with not necessarily analytic coefficients (and as such only
defined on a ray) having possibly divergent series at infinity
(divergent no faster than factorially) and first order polynomially
nonlinear equations. Removing the analyticity and linearity
assumptions raises difficulties requiring new techinques for the
proofs. In particular we provide a method of estimating
the growth of the coefficients of  formal series solutions
based on the recurrence relation that they satisfy.
We mention in addition the paper \cite{Olver} in which very
interesting re-expansions for the optimal remainder are 
obtained in the analytic case.

We also discuss the possibility of extending  our
results to more general systems of differential equations.

 We postpone further  discussions and examples until after Theorem
{\ref{T1}} below.  In Section 2 we study, in some
generality, the connection between optimal asymptoticity and the
method of truncation near the least term. Section 3 is devoted to the
proofs and in Section 4 we discuss, in the context of a particular
family of differential equations, the question of typicality of the
divergence of the asymptotic representations and give a decomposition
formula which provides a different perspective on complete asymptotic
expansions.

Consider the second order linear differential equation

\be{dif2ord}
y''+a(x)y'+b(x)y=0 
\ee

\z given on the ray 

$${\cal R}_{\theta}=\left\{x:e^{-i\theta}x>x_0\right\}$$

The coefficients $a(x)$ and $b(x)$ are assumed to be in $L_{loc}^1({\cal R})$ 
nd to have asymptotic
power series at infinity of the form

\be{asy a,b}
a(x)\sim\sum_{k=0}^{\infty}\frac{a_k}{x^k};\ \ 
b(x)\sim\sum_{k=0}^{\infty}\frac{b_k}{x^k}\ \
(x\,e^{-i\theta}\rightarrow +\infty)
\ee

Moreover, we require that the functions $a(x)$ and $b(x)$ satisfy
a Gevrey-like condition 
\cite{Ramis} on the given ray namely, for some  $\kappa
<|a_0^2-4b_0|^{-\frac{1}{2}}$ 

\ba{Gevrey2}
{a_n,b_n}&<&\mbox{const}\kappa^n n! \cr
|a(x)-\sum_{k=0}^{n-1}\frac{a_k}{x^k}|&<&
 \mbox{const}\kappa^{n}n!\,{|x|^{-n}}\cr
|b(x)-\sum_{k=0}^{n-1}\frac{b_k}{x^k}|&<&\mbox{const}{\kappa^{n}n!}\,{|x|^{-n}}
\end{eqnarray}
{\em uniformly} in $x\in{\cal R}_{\theta} ,n\in \NN$. 

If the conditions (\ref{Gevrey2}) hold, then the rate of divergence
of a formal series solution depends only on the first few terms
in the asymptotic series of $a(x)$ and $b(x)$. The results below
(not merely the proofs) depend on this assumption.

 The polynomial 

\be{poly}
\lambda ^2+a_0\,\lambda +b_0
\ee

\z is assumed to have distinct roots, $\lambda_1, \lambda_2$
(if the roots coincide then there are no terms beyond all 
orders to worry about).
The following expression is a   {\it formal}  solution:

\be{compasy2}
\tilde S=C_1\,e^{\lambda_1\,x }\,x^{r_1}\,\tilde S_1(x)+
C_2\,e^{\lambda_2\,x}\,x^{r_2}\,\tilde S_2(x)
\ee

\z where for $i=1,2$
$$r_{i}=-\frac{a_1\,\lambda_{i}+b_1}{a_0+2\lambda_{i}}$$
 and 
$$
\tilde S_{i}(x)=\sum_{k=0}^{\infty}\frac{s_{i;k}}{x^k}
$$
\z  are formal power series.
The asymptotic behavior for large $n$ of the coefficients of the power
series, as follows from the Proposition {\ref{CvL}} below, has  the form

\be{asymptcoeff}
s_{1;n}\sim\left(R_1+\sum_{m=1}^{\infty}\frac{R_{1,m}}{n^m}\right)
\frac{\Gamma(n+r_2-r_1)}{(\lambda_1-\lambda_2)^n}
\ee

\z (the expression for $s_{2,n}$ is obtained from the one above
by interchanging the indices $1$ and $2$).  Choosing

$$s_{1;0}=s_{2;0}=1$$
the coefficients of the series are uniquely determined.

The expression (\ref{compasy2}) containing  two arbitrary constants
is the general formal solution of our equation, in the  
differential algebra generated by power series and exponentials
(see \cite{Cope}, \cite{Ritt}, \cite{Ecalle}).

The notion of optimal asymptoticity can be extended in a natural way
to asymptotic structures of the form (\ref{compasy2}). We say that
$f$ is uniformly asymptotic to $\tilde S$ for a weight $(w_1,w_2)$
along the ray ${\cal R}_{\theta}$ iff

\be{mainthm}
\left|y(x)-C_1\,x^{r_1}\,e^{\lambda_1\,x}\sum_{j=0}^{k_1-1}\frac{s_{1;j}}{x^j}
-C_2\,x^{r_2}\,e^{\lambda_2\,x}
\sum_{j=0}^{k_2-1}\frac{s_{2;j}}{x^j}\right|<
\ee
\vskip -0.3cm
$$|x^{r_1}e^{\lambda_1\,x}|w_1(k_1)+|x^{r_2}e^{\lambda_2\,x}|w_2(k_2)
\ \ \forall x\in {\cal R}_{\theta} \ \mbox{and }k_1,k_2\in \NN$$

\z We say that $y$ is optimally asymptotic to $\tilde S$ 
along the ray ${\cal R}_{\theta}$ and we write

\be{completeexpan}
y\simeq C_1\,x^{r_1}\,e^{\lambda_1\,x}\sum_{j=0}^{\infty}\frac{s_{1;j}}{x^j}
+C_2\,x^{r_2}\,e^{\lambda_2\,x}
\sum_{j=0}^{\infty}\frac{s_{2;j}}{x^j}
\ee

\z if there is a weight $(w_1(y;n),w_2(y;n))$ such that 
(\ref{mainthm}) holds and moreover for any  $f$ $\in{\cal F}$,
$w_1(y;n)\le w_1(f;n)$ and $w_2(y;n)\le w_2(f;n)$ for all large $n$.

Theorem {\ref{T1}} below states that for each formal
solution there is, on the given ray, a unique true
solution of the equation which is optimally asymptotic to it 
relative to the space of solutions of the same differential equation.
 The theorem can be easily extended to
encompass differential equations that can be brought to the form that
is treated here by an algebraic change of independent variable and
linear changes of the dependent variable; see also Remark 4 below.

Let $B$ be  a (large enough, positive) constant, $r:=r_2-r_1$,
$\rho=\Re(r)$, and $\Delta=\lambda_2-\lambda_1$.

\bt{T1} For 
any $C_1,C_2$ there is
a {\rm unique} true solution $y(x)$ of the differential equation
(\ref{dif2ord}), under the assumptions following it, with the property
(\ref{mainthm}) for the weight

$$w_i(y;k_i)=|C_i|(A_i\sqrt{k_i}+B)\Gamma(k_i-(-1)^i\rho)|\Delta|^{-k_i}\
\ i=1,2$$
\vskip -0.5cm

\z with

\be{defA12}
A_1=\left\{\begin{array}{ll}
a_*|R_1|&\ \ 
\mbox{if $e^{i\theta}(\lambda_1-\lambda_2)\in\RR^+$}\cr
0&\ \ \mbox{otherwise}
\end{array}\right.
\ee

\z where 
 \be{def{a_*}}
a_*:=\sqrt{2}\max_{x\ge  0}e^{-x^2}\int_0^{x}e^{t^2}d\,t=0.765151..
\ee

\z and $A_2$ is  obtained  from (\ref{defA12}) by interchanging
the indices $1$ and $2$.

\end{Th}

\vskip 0.3cm

\vskip 0.3cm
\z {\em Comments.} 

1.)  Let $y(x)$ be the unique solution provided by the theorem and
choose the special truncation orders 
$k_1=k_2=[|x\Delta|]$, which correspond generically to ``truncation to
the least term'' of the series. It will follow from the proof that the
remainder (lhs of (\ref{mainthm})) is in this case less than
\be{truncleastterm}
\mbox{const}\, \frac{1}{\sqrt{|x|}}e^{-|x\Delta|}\max\left\{|x^{r_2}e^{\lambda_1 x}|,
|x^{r_1}e^{\lambda_2 x}|\right\}
\ee

\z i.e., asymptotically smaller than any
(nonzero) solution of the equation by a factor of $\sqrt{|x|}$ on the
Stokes line $\Im (x\Delta)= 0$ and by an exponential factor in the
generic case $\Im (x\Delta)\ne 0$, indicating why there is uniqueness
of the representation.  There is, in effect, a connection between the
shape of the optimal weight and the precision of the least term
truncation method, as will become clear in the folowing section.
Moreover, the
solution that is best approximated by the simultaneous least term
truncation of the two component series in (\ref{completeexpan}) is
precisely the one having the optimal weight. There are differences
between the two approaches but the terms beyond all orders that they
predict agree for this class of differential equations.

\vskip 0.3cm

2.) The values of the constants $A_1$ and $A_2$ in ii) are crucial for
the result.  If say, $\Re(e^{i\theta}\lambda_2)<\Re(e^{i\theta}
\lambda_1)$ and
 $R_1\ne 0$ then there will be no solutions satisfying the required
inequalities with $A_1<a_*|R_1|$ whereas if
$A_1>a_*|R_1|$ there will be infinitely many such
solutions. In contrast, the theorem is true for any large enough $B$.
We also want to stress that uniqueness is {\em relative} to
the space of solutions of the differential equations treated here.

3.) As a manifestation of the Stokes phenomenon (\cite{Stokes},
\cite{Sibuya}, \cite{Wasow}, \cite{Berry}), the constants $C_i$
in the asymptotic representation of a given
 solution depends on the direction $\theta$ of 
the ray considered (provided, of course, that
the differential equation satisfies our hypothesis on more than one
ray at infinity). The example (\ref{eqEi}) given at the begining 
of the section illustrates this point.

The derivation of (\ref{eqEi}) is given at the end of Section 2.  We
also mention at this point that a more detailed analysis is possible
along the same lines and it shows that the term beyond all orders
varies smoothly on a scale of the order $\arg(x)\sim x^{-1/2}$ in a
way which agrees with the results obtained with the hyperasymptotic
technique of Berry \cite{Berry}, \cite{Berry-hyp}. We will however not
pursue this issue here.

4.) One can in principle allow for more general formal structures than
(\ref{compasy2}), for instance those obtained by substitutions and
formal algebraic operations on (\ref{compasy2}) (for a discussion on
what formal structures are relevant to solving differential equations
see \cite{Ecalle}), to allow for asymptotic representations of the solutions
of equations that are not of the form required by the theorem, but can
be brought to that form; the asymptotic inequalities are well suited
for simple algebraic operations. Consider for instance
homogeneous Airy equation

\be{Airyomo}
y''-x\,y=0
\ee

\z  After the substitution
$y(x)=\exp(\frac{2}{3}\,x^{\frac{3}{2}})g(x)$ 
followed
by taking $x=s^{2/3}$ we get

\be{Aiges}
g''+(\frac{4}{3}+\frac{1}{3s})g'+\frac{2}{9s}\,g=0
\ee

\z to which Theorem \ref{T1} applies and we obtain, after undoing
the transformations, the
following asymptotic representation of the general solution of the
Airy equation (\ref{Airyomo}):

\be{AsyAirygen}
y(x)\simeq C_1e^{\frac{2}{3}x^{\frac{3}{2}}}x^{-1/4}
\sum_{k=0}^{\infty}\frac{s_k}{x^{\frac{3k}{2}}}+
 C_2e^{-\frac{2}{3}x^{\frac{3}{2}}}x^{-1/4}
\sum_{k=0}^{\infty}\frac{(-)^k s_k}{x^{\frac{3k}{2}}}
\ee

\z when $x$ becomes large along a given ray in the complex plane. For a
fixed solution $C_1,\ C_2$ will depend on the ray. A brief
derivation of (\ref{AsyAirygen}) and the expressions of $A_1,\ A_2$
and $s_k$ are given at the end of Section 2.

5.) Whenever
$e^{i\theta}\Delta\notin\RR^-$ all the terms of the
expansion (\ref{compasy2}) are simultaneously visible to the
inequalities, as it will follow from the proof
of the theorem.  If on the other hand
$e^{i\theta}\Delta<0$ i.e., on the Stokes line, besides
the first series only the first term of the second series is caught by
the inequalities. However this first term is enough in order to
establish a one-to-one (linear) correspondence between formal and true
solutions, so that we keep, by courtesy, all the other terms of the
second series.

6.) For obtaining  complete asymptotic representations
by means of optimal inequalities  for higher order differential equations
or systems of equations it seems that
it is necessary to use, inductively, asymptotic comparisons with solutions
of lower order ODE's.

\vskip 0.3cm

The result below shows that a given function cannot be a solution
of two essentially different second order differential equations
that satisfy our assumptions.
We say that a differential equation of the type considered
in the theorem 
is in the canonical form if:

\be{simpli}
\theta=0;\ \lambda_1=0;\ |a_0|=1;\ r_1=0;\ b_0=b_1=0
\ee

This can be always achieved by the substitution $y(x)=\exp(\lambda_1\,x)x^{r_1}\tilde y(x)$
followed by a change of independent variable 
$x=e^{i\theta}|\lambda_1-\lambda_2|^{-1}\tilde x$. The differential
equation is said to be formal if the coefficients $a(x),b(x)$ are
replaced by their formal series.

The condition (\ref{Gevrey2}) plays a crucial role in the result below.

\bp{Puniq} Let $\tilde S=\sum_{n=0}^{\infty}s_n
x^{-n}$ be a formal power series such that $s_n$ have the asymptotic
behavior (\ref{asymptcoeff}) with $R_1\ne 0$. 
Then $\tilde S$ formally solves at most one second-order
{\em canonical} formal differential equation of the form
(\ref{dif2ord}) if the formal series of $a(x)$ and $ b(x)$ satisfy
the first condition in (\ref{Gevrey2}).
\end{Pro}

It follows for instance that within the class of differential
equations we are concerned with, there is a one-to-one correspondence
between generic formal solutions and true solutions.  Indeed, a true
solution is associated with a unique representation within the
solutions of the differential equation it originates in.  But then,
the representation determines uniquely the power series of the
coefficients of the equation and by this determines uniquely the
differential equation itself. Conversely, a formal expansion {\em
solution} is associated with a unique differential equation and within
that differential equation with a unique true solution.  In this sense
we also obtain a summation method.  One can attempt to continue this
construction inductively, considering at the $n$-th step differential
equations whose coefficients are solutions of the equations gotten at
the step $n-1$ and obtain a more general representation spaces and
summation methods. 
\vskip 1ex 
{\centerline{*}}
\vskip 1ex 

\z Similar results hold for first order nonlinear differential
equations of the form

\be{1stnonln}
y'=F(y,x):=\sum_{p=0}^{P}{F_p(x)}{y^p}
\ee

\z defined on the ray ${\cal R}_{\theta}$.

We require that for  $p=0,..,P$, $F_p(\cdot)\in C[R_{\theta}]$,
$F(0,\infty)\equiv
f_{0,0}=0,\ 
F_y(0,\infty)=:\alpha=-\ f_{1,0}\ne
0$ and impose the Gevrey-like condition

\ba{ineqf_i}
\left|F_p(x)-\sum_{k=0}^{n-1}\frac{f_{p,k}}{x^k}\right|&<&
               C_F{\kappa^{n}n!}{|x|^{-n}}\cr
|f_{p,n}|&<&C_F\kappa^n\,n!
\end{eqnarray}

\z (where $C_F$ is a constant) uniformly in $x\in{\cal R}_{\theta}, n\in\NN$, where $\kappa<|\alpha|^{-1}$.

We are only interested in the case when the
equation allows for exponentially {\em small} terms beyond all orders;
the condition for that (as it will become immediately clear) is

\be{cond on f_{1,0}}
\arg(x\alpha)\in (-\frac{\pi}{2},\frac{\pi}{2})
\ee

\bl{lemma5}
There exists a formal series solution 
\be{defseri}
S:=\sum_{k=1}^{\infty}\frac{s_k}{x^k}
\ee

\z of the equation (\ref{1stnonln}).
The coefficients of the series have the asymptotic 
behavior 
\be{asybeha}
s_k\sim\left(R+\sum_{m=1}^{\infty}\frac{R_m}{k^m}\right)\alpha^{-k}\Gamma(k+r)
\ \ ({\mathrm as\ }k\rightarrow\infty)
\ee

\z where 

$$r=f_{1,1}+2\,f_{2,0}f_{0,1}/\alpha$$
\end{Le}

\z If we are now looking for nearby formal solutions 
$S_{\delta}:=S+\delta$ we get 

\be{lineariz. asym.}
\delta'=\left(-\alpha+\frac{r+1}{x}\right)\delta  +O(\delta^2,\delta/x^2)
\ee

\z so that for large $|x|$

\be{asym. delta}
\delta\sim C\,x^{r}e^{-\alpha x}
\ee

The analog of Theorem \ref{T1} in this case reads: consider the
equation (\ref{1stnonln}) under the assumptions
(\ref{1stnonln})-(\ref{cond on f_{1,0}}) and let
$\rho=\Re(r)$. Let $B$ be large enough. 

There is a one-to-one asymptotic correspondence between
true solutions which decay at infinity and formal
solutions:

\bt{T2} Given any constant $C$   there
exists a {\rm unique} true solution $y(x)$ of the differential equation
(\ref{1stnonln}) such that 
\be{eqT2}
\left|y(x)-\,\sum_{j=0}^{k-1}\frac{s_{j}}{x^j}
 -C\,x^{r}\,e^{-\alpha\,x}
\right|<(A\sqrt{k}+B)\Gamma(k+\rho)||x|^{-k} 
\end{equation}

\z for all $x\in {\cal R}_{\theta}$ and  $k\in \NN$, where 
\be{defA}
A=\left\{\begin{array}{ll}
a_*|R|&\ \ 
\mbox{if $e^{i\theta}\alpha\in\RR^+$}\cr
0&\ \ \mbox{otherwise}
\end{array}\right.
\ee

Conversely, given a solution of (\ref{1stnonln}) such
that $xy(x)$ is bounded for large $x$, there is a 
(unique) constant $C$ such that, for all $x\in {\cal R}_{\theta}$ and
$k\in \NN$, (\ref{eqT2}) holds.

\end{Th}

\z In the sense of (\ref{eqT2}) we then write

\be{rep1st}
y(x)\simeq \sum_{j=0}^{\infty}\frac{f_{j}}{x^j}+
C\,x^{r}\,e^{-\alpha\,x}
\ee

The comments that follow Theorem {\ref{T1}}, with
obvious adaptations,  also apply in this case.

We could actually continue  the construction of the formal solution
and consider the complete formal solutions or
``transseries''(\cite{Ecalle}
,\cite{Ecalle2})
which in this case have the form

$$Y_0+x^{r}e^{-\alpha\,x}Y_1+x^{r_2}e^{-2\,\alpha\,x}Y_2+...$$

\z in which $Y_i$ are formal power series. $Y_1$ is determined up
to an arbitrary constant which, once given, determines completely all
the following power series $Y_i$, $i\ge 2$. We will however not pursue
this direction here but merely remark that the one-to-one
correspondence between formal solutions and true solutions is pinned
down by the first term of the second formal series and contend to
control the asymptotics only to that level.

\end{section}

\begin{section} {Optimal estimates for a class of divergent series}

In this section we estimate
the optimal bounds  for series that diverge in a way typical
for differential equations and make the connection between
optimal uniform inequalities and the method of truncation
to the least term. We also give results on  the precision 
with which a function is represented by its optimal
asymptotic expansion. We consider formal power series, for large argument
on a ray and establish the connection to the technique 
of optimal truncation of divergent series.

\be{psG}
\tilde S:=\sum_{k=0}^{\infty}\frac{c_k}{x^k}
\ee

\z where 
\be{defgenps}
c_k=(1+\epsilon_k)\,k^{i\phi}C^k\exp(G(k));\ \ \epsilon_k\rightarrow 0
\ \ as\ \ k\rightarrow\infty
\ee

\z Since the term $C^k$ can be always absorbed into the independent variable
so that we will take  $C=1$.

In (\ref{defgenps})  $\phi$ is real and we make the following
assumptions.  $G(\cdot)\in C^3(\RR^+)$ is a real-valued increasing,
convex function with the properties $G''(x)\rightarrow 0$ and
$x^2G'''(x)\rightarrow -\gamma<0$ as $x\rightarrow\infty$.  
Such a series is necessarily divergent for all values
of $x$ (convergent series are discussed in the introduction).
An example satifying the requirements would be  
$c_k=C^k\Gamma^{\beta}(k+r)$.

The behavior of the optimal weight for a divergent series
depends critically on the ray along which we consider
the asymptotic series. In the setting (\ref{defgenps})
the optimal weight is larger by a factor
of order $\sqrt{n}$ along the real positive axis than it is
along any other ray. This is a manifestation
of the Stokes phenomenon at the level of the asymptotic
series themselves. Intuitively we can account for this behavior in the
following way. A natural scale for studying the difference
between a function and the $n-th$ truncate of its asymptotic series
is the $n+1$-th term of the series (indeed, the $n+1$-th term is meant to be a
correction for the above mentioned difference). When $x$ is very
large, the successive terms of the series start by decreasing fast
as so does the error in approximating the function by its series.
The least error suggested by our rough guide, the next term
of the series, reaches a minimum when

\be{eq.l.t.}
G'(s_x)=\ln(|x|)
\ee

The width of the minimum is of the order $\sqrt{s_x}$ and within this
width, for $x>0$, the ratio of two terms is approximately one.  But
then we realize that if at one point within this range the difference
between the function and the truncate of the series is approximately
equal to the next
term there will be a point within the same region where the difference
will be roughly $\sqrt{n}$ times as big. If $x$ is not real the ratio
of the successive terms is of the form $e^{-i\theta}$ so that the
overall accumulation of errors is still of the order of a constant. 

 We want to stress that, as it follows from Lemma \ref{lemmaex} below
and the estimates in its proof, that  the least-term truncation of a
formal series solution is a good approximation only for {\em one}
solution of a given differential equation; for all all the other
solutions there are exponential corrections that are much larger than
the least term of the series and have to be removed before calculating
the function from its series (to such an accuracy).

Consider the power series $\tilde S$, under the given assumptions,
for $x$ on the ray ${\cal R}_{\theta}$.

\bl{OptimSer}
i) For any $\theta$ there exist (smooth) functions that
are optimally asymptotic to the series (\ref{defgenps}). 
A optimal weight is
$$g_*(n)=a_*\gamma^{-\frac{1}{2}}(\sqrt{n}+B){|c_{n+1}|}\ 
{\mbox{for }} \theta=0$$
\z More precisely, assume that $\Phi$ is a function
with the property that $\exists \eta $ and  $B$ such that 
$\forall\ x\in {\cal R}_{\theta}$ and $\forall n\in \NN$ we have

\be{fundamineq}
\left|\Phi(x)-\sum_{j=0}^{n}\frac{c_j}{x^j}\right|
<\eta \,a_*\gamma^{-\frac{1}{2}}(\sqrt{n}+B)\frac{|c_{n+1}|}{x^{n+1}}
\ee
Then, $\eta \ge 1$ and there exist functions satisfying
(\ref{fundamineq}) for $\eta=1$, for some $B$.

ii) For $\theta \ne 0$ assume further that

\be{cond-eps}
\epsilon_k=o(k^{-1/2})
\ee

\z Then the optimal 
weights can be taken
$$g_*(n)=(a(\theta)+Bn^{-1})|c_{n+1}|$$

\z where the constant $a(\theta)$ can be computed
explicitly and for small $\theta$ has the behavior
\be{1/theta}
a(\theta)\sim \theta^{-1}
\ee

\z The behavior (\ref{1/theta}) is actually present
in the concrete estimates of the exponential integral, see 
(\ref{estiEi,away}).

\end{Le}
\vskip 0.5cm

\z 
The result below answers the
existence part of Lemma {\ref{OptimSer}}, i) and gives the connection
with the technique of summation to the least term. Let $n_x$  
be defined, for
$x$ large enough by $n_x=[s_x]$, (cf. (\ref{eq.l.t.});\ [\ ]
denotes the integer part)

\bp{OptSummationProp}

\z i) Let $\theta=0$. A function $\Phi$ satisfies the uniform inequalities
(\ref{fundamineq}) with $\eta =1$ iff

\be{OptSumProp1}
\limsup_{x\rightarrow\infty}
n_x^{-1/2}\frac{x^{n_x}}{|c_{n_x}|}
\left|\Phi(x)-\sum_{j=0}^{n_x}\frac{c_j}{x^j}\right|=0
\ee

\z Instead, if the limit above is $\epsilon>0$ and
$\Phi$ satisfies (\ref{fundamineq}) for some $\eta $ then $\eta
\ge 1+\epsilon$.

ii) If $\theta\ne 0$ and (\ref{cond-eps}) holds, then a function
$\Phi$ is uniformly  asymptotic to the series (\ref{defgenps}) along
the ray ${\cal R}_{\theta}$ with respect to a weight

$$w(n)=C|c_{n+1}|$$
 for some $C$
iff

\be{OptSumProp2}
\limsup_{|x|\rightarrow\infty}
\left|\frac{x^{n_x}}{c_{n_x}}\right|
\left|\Phi(x)-\sum_{j=0}^{n_x}\frac{c_j}{x^j}\right|<\infty
\ee

\end{Pro}

\vskip 0.5cm

\z {\em{Corollary}}. If $\Phi_1$, $\Phi_2$ are two functions
satisfying the hypothesis of part i) of $\mbox{Lemma \ref{OptimSer}}$
with $\eta=1$ then

$$
\limsup_{x\rightarrow\infty}
n_x^{-1/2}\frac{x^{n_x}}{|c_{n_x}|}
|\Phi_1(x)-\Phi_2(x)|=0
$$

\z In particular when $G(n)=\ln(\Gamma(n+n_0))$ then

\be{uniq;coro}
\Phi_1(x)-\Phi_2(x)=o(e^{-x}\,x^{n_0})
\ee

\z Correspondingly, in the case iii) we have 

$$
\limsup_{x\rightarrow\infty}
n_x^{-1/2}\frac{x^{n_x}}{|c_{n_x}|}
|\Phi_1(x)-\Phi_2(x)|=0
$$

\z and for $G(n)=\Gamma(n+n_0)$,

\be{uniq;coro;compl}
\Phi_1(x)-\Phi_2(x)=O\left(e^{-|x|}\,x^{n_0-1/2}\right)
\ee

Classically, two functions $\Phi_1$ and $\Phi_2$  have the
same asymptotic series (provided they have one) if their difference is
asymptotically less than any power of $x$. The corolarry shows that
there is a definite precision gain when requiring optimal uniform
estimates.

\end{section}

\begin{section}{Proofs and further results} 

Before giving the general proof we mention as an illustration the
particularly easy computation of the optimal weight for the
exponential integral: see the arguments
starting with eq. (\ref{eig(n)start}).

In this section
we make the following convention: for $n>k$,

$$\sum_{j=n}^k a_j=-\sum_{j=k}^n a_j$$

For the proof of Lemma \ref{OptimSer} we need the following
elementary result. The notations and hypothesis are those
that precede  Proposition
\ref{OptSummationProp}.
\bp{uniform,factorial} 
i)Let  $k_0$ be large enough.  Then for $\theta=0$,

\be{supest}
\lim_{x\rightarrow\infty}\,\,\sup_{K\in \NN}
\left\{(K+k_0)^{-1/2}{{x}^{K}}{e^{-G(K)}}\left|
\sum_{j=n_x}^{K}\frac{c_j}{x^{j}}\right|\right\}=a_*\gamma^{-\frac{1}{2}}
\ee

\z The supremum in  (\ref{supest}) is actually attained
for
 
$$K\sim n_{\pm}:= [s_x\pm\frac{1}{a_*}\sqrt{\frac{\gamma}{s_x}}]$$
\vskip 0.3cm
\z ii) If $\theta\ne 0$ then

\be{supest,complex}
\lim_{x\rightarrow\infty}\,\,\sup_{K\in \NN}
\left\{{{x}^{K}}{e^{-G(K)}}\left|
\sum_{j=n_x}^{K}\frac{c_j}{x^{j}}\right|\right\}<\infty
\ee

\end{Pro}

{\em Proof

} The proof is essentially straightforward. 
Note first that from the assumptions
it follows that $G''(x)\sim \gamma/x$ for large $x$.

Note also that we need only  consider the case $\epsilon_j=0$. Indeed,
the presence of a (large enough) constant $k_0$ makes the 
proposition above insensitive to the behavior
of $c_j$ for small $j$. On the other hand,

\be{elim:epsilon}
\left|\sum_{j=n_x}^{K}\frac{c_j}{x^j}\right|=
\left|\sum_{j=n_x}^{K}j^{i\phi}e^{-G(j)+j\,\ln(x)}\right|+
O(\max_{K\le j\le n_x}|\epsilon_j|))\sum_{j=n_x}^{K}e^{-G(j)+j\,\ln(|x|)}
\ee

\z and the maximum above approaches zero as $K\rightarrow\infty$.

$\bullet$ i) Take first $\phi=0$. Let $\beta\in(\frac{1}{2},\frac{2}{3})$ and
start with the range of $K$ so that $$|K-n_x|<n_x^{\beta}$$ which is
actually the important range with respect to (\ref{supest})

Here,
the Euler-Maclaurin summation method is suited. For definiteness
we take  $K\ge n_x$ (the other case is very similar). We have

\ba{eumacfi0}
&&\frac{x^{n_x}}{c_{n_x}}e^{G(j)-j\ln\, x}=
(1+O(s_x^{3\beta-2}))e^{\frac{1}{2}G''(s_x)(j-s_x)^2}=\cr
&&(1+O(s_x^{3\beta-2}))\int_{j}^{j+1}e^{\frac{1}{2}G''(s_x)(j-s_x)^2}dt=\cr
&&(1+O(s_x^{3\beta-2}))\int_{j}^{j+1}e^{\frac{1}{2}G''(s_x)(t-s_x)^2}dt\cr
&&
\end{eqnarray}

\z so that

\ba{finev}
&&\hskip -10pt\max_{|K-n_x|<n_x^{\beta}}\left\{K^{-1/2}e^{-G(K)+K\,\ln\, x}
\sum_{j=n_x}^{K}e^{G(j)-j\,\ln\, x}\right\}=\cr
&&\hskip -10pt\left(1+O(n_x^{3\beta-2})\right)\max_{|K-n_x|<n_x^{\beta}}
\left\{K^{-1/2}e^{-\frac{1}{2}G''(s_x)(K-s_x)^2}
\int_{n_x}^{K+1} e^{\frac{1}{2}G''(s_x)(t-s_x)^2}d\,t\right\}=\cr
&&\hskip -10pt\left(1+O(n_x^{3\beta-2})\right)a_*({K\,G''(s_x)})^{-\frac{1}{2}}
\rightarrow a_*\gamma^{-\frac{1}{2}}
\end{eqnarray}

\z as $x\rightarrow\infty$.

$\bullet$ It remains to obtain an upper bound in the region
$|K-n_x|\ge n_x^{\beta}$. Assume $K<n_x$ (the opposite
case is treated similarly) and let 
$q=s_x-\frac{1}{2}s_x^{\beta}$. With $F(z):=G(z)-z\,\ln(x)$ we have in
view of the estimates above,

\be{splitsum}
\sum_{j=K}^{n_x}e^{F(j)-F(K)}
\le
\sum_{j=K}^{q}e^{F(j)-F(K)}+\mbox{const}\sqrt{q}
e^{F(q)-F(K)}
\ee

\z On the other hand, from the definition of $s_x$,

$$F(j)-F(K)=\int_{K}^{j}dz\int_{s_x}^{z}d\,t\,G''(t)\le\,
\mbox{const}\int_{K}^{j}\ln(z/s_x)dz$$

\z Throughout the domain of integration $s_x>z+s_x^{\beta}>z+z^{\beta}$
so that

$$\ln(z/s_x)<\ln\left(1-\frac{z^{\beta}}{z+z^{\beta}}\right)
<-\frac{1}{2}z^{\beta-1}$$
and we get, for some positive constants $C_1,C_2$,
$$
\sum_{K}^{q}e^{F(j)-F(K)}<
\sum_{n=0}^{q-K}e^{-C_1[(K+n)^{\beta}-K^{\beta}]}<
C_2\,(K+k_0)^{1-\beta}<
\frac{1}{2}\sqrt{\frac{(K+k_0)}{\gamma}}a_*
$$

\z In the same way, for some positive constants,

\ba{final i)}& e^{F(q)-F(K)}\le e^{-C_1[q^{\beta}-K^{\beta}]}<
e^{-C_2 K^{\beta-1}(q-K)}<e^{-C_3 q^{\beta}}
<e^{-C_4 x^{\beta/\gamma}}\cr
&
\end{eqnarray}

\z so that the second term in (\ref{splitsum}) vanishes as 
$x\rightarrow\infty$.

$\bullet$ If now $\phi\ne 0$ the upper bounds are trivially obtained by taking
absolute values in the sums. For a lower bound we note simply that

$$\sum_{n_x}^{n_+}\frac{c_j}{x^j}=n_x^{i\phi}\sum_{j=n_x}^{K}
e^{-G(j)+j\,\ln(x)}\left(1+O(n_x^{-1/2})\right)$$

$\bullet$ ii) Let again $F(k)=G(k)-k\ln(|x|)-i k\theta$. 
Given a function
$f(x)$ 
which is asymptotic to the series $\tilde S$ we write

\be{decomp f}
f(x)=:\sum_{k=0}^{n_x}\frac{c_k}{x^k}-
\frac{e^{F(n_x)}}{1-e^{-i\theta}}\chi(x)
\ee

\z and we will choose $\chi$ to make $f(x)$ optimally asymptotic to 
$\tilde S$.
Consider first the region $|k-n_x|<n_x^{\beta}$ and define

\be{equF}
\sigma_k:=\frac{k^{i\phi}e^{F(k)}}{(e^{-i\theta}-1)}
\left(1-\frac{\gamma e^{i(k-n_x)}}{n_x}
\sum_{j=-n_x^{\beta}}^{K-n_x-1}je^{-ij\theta}\right)
\ee

\z It is easy to check that

\be{guessEuMac}
\sigma_{k+1}-\sigma_k=\frac{c_k}{x^k}\left(1+o(n_x^{-1/2})\right)
\ee

\z (to get (\ref{equF}), it is simpler in this case to solve
(\ref{guessEuMac}) directly by perturbation expansion than to use the
Euler-Maclaurin summation formula). Using (\ref{decomp f}),
(\ref{equF}) and (\ref{guessEuMac}) we get

\ba{formerr}
E_n:=\left|\frac{x^{n+1}}{c_{n+1}}\right||f(x)-\sum_{k=0}^n\frac{c_k}{x^k}|=
\qquad\cr
\frac{1}{|1-e^{-i\theta}|}
\left|1-e^{i (n_x-n)\theta}
\left(e^{-i\theta h_x(n)}-\chi+o(1)\right)
e^{-\gamma\frac{(n_x-n)^2}{2n_x}}\right|
\end{eqnarray}

\z where $h_x(n)=1$ if $n>n_x$ and is zero otherwise.
For large $x$ (\ref{formerr}) can be visualised as distance between
$z=1$ to the points on two spirals centered at $e^{-i\theta
h_x(n)}-\chi$, and with diameter slowly decresaing as $|n-n_x|$
increases.  It is not hard to see that there is in this region a best
choice $\chi_*$ of $\chi(x)$ (which minimizes the maximal error in
(\ref{formerr}) since for any value of $x$ there is only a finite set
of $n$ such that $|n-n_x|<n_x^{\beta}$. For $\theta$ incommensurate
with $\pi$ the geometry of the problem shows that

\be{chi*}
\chi_*=e^{-i\frac{\theta}{2}}\cos\frac{\theta}{2}+o(1)
\ee

\z for which we get

\be{a(theta)}
a(\theta)=\frac{1+|\sin\frac{\theta}{2}|}{|1-e^{-i\theta}|}=
\frac{1}{2}\left(1+\frac{1}{|\sin\frac{\theta}{2}|}\right)
\ee

\z If $\theta$ is a rational multiple of $\pi$ the computation 
can still be done explicitly, in a straightforward manner and 
$a(\theta)$ is slightly less than (\ref{a(theta)}).

In the region $|k-n_x|>n_x^{\beta}$ a similar calculation shows that
the error is bounded by

$$E_n \le |1-e^{F'(n)-i\theta}|^{-1}+o(1)|\chi(x)|\le a(\theta)+o(1)|\chi(x)|$$

\z which shows that $\chi_*$ is indeed the optimal global choice
of $\chi$, as it is not possible to decrease the error in this region
without making it unbounded in the $n-n_x|<n_x^{\beta}$ region.

$\square$

{\em{Proof}} of Lemma \ref{OptimSer}. For part i), if we assume there
existed a function $\Phi$ satisfying the estimate (\ref{fundamineq})
uniformly in $x$ and $n$ for some $\eta <1$ it would follow that

\be{eq1000}
\pm\left(\Phi(x)-\sum_{j=0}^{n_{\pm}}\frac{c_j}{x^j}\right)\le 
\eta \,a_*\gamma^{-\frac{1}{2}}(\sqrt{n_{\pm}}+B)\frac{c_{n_{\pm}+1}}{x^{n_{\pm}+1}}
\ee

\z By adding the two inequalities above we arrive at 
an immediate contradiction with Proposition {\ref{uniform,factorial}}.

Part ii) and iii)  are straightforward consequences of Proposition
{\ref{OptSummationProp}}. Indeed,  note first
an easy example 
of a  function satisfying (\ref{OptSumProp1}) and (\ref{OptSumProp2}):

\be{def f. by s.l.t}
\Phi_*(x):=\sum_{j=0}^{n_x}\frac{c_j}{x^j}
\ee

It is not dificult to smooth $\Phi_*$ and preserve
(\ref{OptSumProp1}), (\ref{OptSumProp2}) since the size of a jump at a
point of discontinuity of $\Phi_*$ is $|c_{n_x}x^{-n_x}|$.

On the other hand, Proposition {\ref{OptSummationProp}}
follows immediately from Proposition \ref{uniform,factorial} by 
the triangle inequality. $\square$
\vskip 1ex
{\centerline{*}}
\vskip 1ex
{\em{Proof}} of Theorem \ref{T1}

There are no assumptions that would distinguish between the quantities
with subscript ``1'' from those with subscript ``2''.  The existence
part of the theorem follows (trivially, using
triangle inequalities) if we prove it for $C_1=1$ and
$C_2=0$. This case will follow from Lemma~{ \ref{lemmaex}} below.

We can assume that the differential
equation is brought to its canonical form because the asymptotic
inequalities are transformed in an obvious way in the substituions
involved.  Also we note the following
inequality that we will use frequently 
and which is clear from the integral
representation of the Gamma function:
$$|\Gamma(x)|\le \Gamma(\Re(x))\ \mbox{for $\Re(x)>-1$}$$

Let $r:=r_2=-a_1,\ \rho:=\Re(r)$ and $ \alpha:=a_0$.

\bl{lemmaex}
There is a solution $Y$ of the differential equation  (\ref{dif2ord})
under the assumptions (\ref{asy a,b})-(\ref{simpli}) which satisfies
the inequalities:
\be{finalineq}
|Y-\sum_{k=0}^{n}\frac{s_{1;k}}{x^k}|<\mbox{const}\, e^{-x}\,x^{\rho-1/2}\ \  for\ \  x\in(n,n+1)
\ee
where {\em const} does not depend on $n,x$
\end{Le}

Before we prove this Lemma we need some results on
the formal solutions.  Let $s_k=s_{1,k}$. The condition
that $\sum_{k=0}^{\infty}\frac{s_k}{x^k}$ is a formal 
solution of the equation leads to the following recurrence relation
for the coefficients of the series:

\be{recp=1}
s_n=\frac{1}{\alpha}\left(n-1+r+\frac{a_1+b_2}{n}\right)s_{n-1}
+\sum_{n\ge j\ge 2}\frac{1}{n\,\alpha}
\left((j-n)a_j+b_{j+1}\right)\,s_{n-j}
\ee

\bp{CvL} The behavior of $s_n$ for large $n$ is 

\be{asymptcoefLin}
s_n\sim\left(R+\sum_{m=1}^{\infty}\frac{R_m}{n^m}
\right)\alpha^{-n}\Gamma(n+r)
\ee

\end{Pro}

In the assumption of analyticity at infinity of
 the coefficients of the equation, the leading behavior
of $s_n$ was shown in \cite{Jurkat}.

{\em{Remark}}. Although
in  a typical case the constant $R$ 
will be  nonzero  as it is easy to understand
by examining the recurrence (\ref{ecred}) below, massive cancellations
are possible so that the series $R+\frac{R_1}{n}+..$
could be zero to all orders in $1/n$;
this case is important in its
own right--see Section 3 for some examples and applications in this
connection.

{\em Proof

}

In order to avoid dealing with the poles of the $\Gamma$
function which might occur for (uninteresting) small values of $n$
we let

\be{defGam_}
\Gamma_{x}:=\left\{\begin{array}{ll}
\Gamma(x)&\mbox{if $\Re(x)>-1$}\cr
1&\ \ \mbox{otherwise}
\end{array}\right.
\ee
\z It is convenient to pull out the leading behavior suggested by (\ref{recp=1}):

$$s_n=:\alpha^{-n}\Gamma_{n+r}c_n$$

\z We get for $c_n$

\be{ecred}
c_n=
c_{n-1}\left(1+\frac{A_n}{n^2}\right)+O(n^{-2})c_{n-2}
+\sum_{n-n_0\ge j\ge 3}C(n,j)c_{n-j}+o({n^{-2}})
\ee

\z where the term with $j=3$ and 
$n_0>\max\{\rho,3\}$ terms with $j$ near $n$  in the sum 
were treated separately. Due to (\ref{Gevrey2}) the latter give a
collective a contribution of $o({n^{-2}})$ and also, the following
estimate holds for $C(n,j)$:

\be{estiC}
|C(n,j)|<\mbox{const}\,\frac{(n-j)!j!}{n!}=O(n^{-3})
\ee
in the range $n-3\ge j\ge 3$ (the inverse of
the binomial coefficient is convex in $j$). 

By induction on $n$ we see that
that $|c_n|<A\,\prod_{j\le n}(1+\mbox{const}\,j^{-2})$ for some
large enough $A$ and $\mbox{const}$, whence the sequence ${c_n}$ 
is bounded. But then it follows immediately from
(\ref{ecred}) that the sequence $\{c_n\}$  is convergent. Furthermore,
taking $c_n=d_n\prod_{j\le n}(1+A\,j^{-2})$ we see that

\be{asdn}
d_n=d_{n-1}+\frac{K_n}{n^2}
\ee

\z where $K_n$ is a bounded sequence and thus,
if $R$ is the limit of the $c_n$ we have 

\be{step1,sn}
c_n=R+O(n^{-1})
\ee

\z It is now easy to bootstrap the estimates for $c_n$ in the recurrence
(taking out explicitely more and more terms from the sum) 
 to get $c_n=R+\frac{R_1}{n}+O(n^{-2})$ and so on.

\z 
$\square$

\z Let

$$U_n=\left(|R|+O(n^{-1})\right)\Gamma(n+\rho)$$

\z represent an upper bound of $|s_n|$. Note that
if $R\ne 0$ then  $|s_n|\sim U_n$ for large $n$.

\z Let 

$$Y_n:=\sum_{k=0}^n\frac{s_k}{x^k}$$

\z  Take $n$ to be large enough
so that for $x>n$ the asymptotic estimates for $a(x),\ b(x)$ 
hold and so that $n+\rho>0$.

\bp{RemainderEq} In the differential equation verified 
by $Y_n$ 
\be{difeqYn}
Y_n''+a(x)Y_n'+b(x)Y_n=R_n
\ee

\z the inhomogeneous term satisfies the estimate

\be{EstRem}
R_n(x)=(R+O(n^{-1}))\alpha^{-n-2}\frac{\Gamma(n+r+2)}{x^{n+2}}
\sim(\sqrt{2\pi}R+O(x^{-1}))\alpha^{-n-2}x^{r-\frac{1}{2}}e^{-x}
\ee
for $x\in(n,n+1)$.
\end{Pro}
{\em Proof

}
We substitute
$a(x)=\sum_{0}^{n}\frac{a_k}{x^k}+\frac{A_{n+1}}{x^{n+1}}$ (and the
corresponding expression for $b(x)$) in the expansion of the LHS of
(\ref{difeqYn}).  To simplify the notations we let throughout this
proof $a_{n+1}=A_{n+1}(x)$; $b_{n+1}=B_{n+1}(x)$. The coefficients of
$x^{-j}$ with $j\le n+1$ vanish by the definition of $Y_n$.  The
coefficients of $x^{-n-1}$ add up to $\alpha (n+1)s_{n+1}$.  Indeed,
$-\alpha (n+1)s_{n+1}$ is the only term of order $n+1$ which is
missing in (\ref{EstRem}) with respect to the corresponding expression
of $Y_M,\ M>n$. So,
\be{defRn}
R_n=
       \alpha\frac{(n+1)s_{n+1}}{x^{n+2}}
+\sum_{m\ge n+3}x^{-m}T_m
\ee
\z with
\be{deftm}
T_m:=
             \sum_{\stackrel{\scriptstyle k+j=m}{k,j\le
n}}\left((k-1)s_{k-1}a_j+s_k\,b_j\right)
\ee
\z An individual term in $T_m$ can be bounded by:

$$\left|(k-1)s_{k-1}a_j+s_k\,b_j\right|\le \mbox{const}\,U_n$$

\z as it is
easy to see using estimates for $s_k,\,a_k,\,b_k$.  Thus
$$x^{-n-3}T_{n+3}<\mbox{const}\, n^{-n-1}U_{n}$$ and the same is
true for $\sum_{m>n+3}x^{-m}T_m$, as shown by
a gross majorization by the number of terms ($O(n)$ by (\ref{deftm}))
times the maximal
term. It follows that the modulus of the sum on the RHS of (\ref{defRn})
is at most 
$$O(n^{-1})\frac{\Gamma(n+r+1)}{x^{n+1}}$$
for $x\in(n,n+1)$. The rightmost estimate in (\ref{EstRem})
is just Stirling's formula.
$\square$

\z We return to Lemma {\ref{lemmaex}}.

{\em Proof

} The equivalent vectorial equation
is now preferable:
\be{eqvec}
{\bf F}'=A\,{\bf F}
\ee
where 
$$A:=\pmatrix{0 & 1 \cr
{-b(x)} & {-a(x)} };\ \ {\bf F}:=\pmatrix{f\cr f'}
$$

\z On the interval $[n,n+1]$ we look for  solutions in the form
$$
{\bf F}=\pmatrix{Y_n\cr Y'_n} +{\bf H}(n,\cdot)
$$

\z ${\bf H}$ has to satisfy the differential equation
\be{dequH}
{\bf H}'=A\,{\bf H}+{\bf T}_n;\ \ {\bf T}_n=\left(
\begin{array}{cccc}
0 \\
-R_n \end{array}
\right)
\ee

\z and the continuity condition (which is actually enough
to ensure the smoothness of ${\bf Y}$)

\be{conti}
{\bf H}(n;n_+)={\bf H}(n-1;n_-)-{\bf E}_n\ \ \ (n_{\pm}:=n\pm 0)
\ee
\vskip -0.3cm
$$
{\bf E}_n:=\frac{s_{n}}{n^{n}}
\pmatrix{1\cr -1}
$$

A differential equation of the form (\ref{dif2ord}) under
the given assumptions admits always two special solutions
which for large $|x|$ have the behavior (see \cite{Wasow})

\be{asyfuset}
\quad\ y_1(x)\sim 1+\frac{s_1}{x}+O(x^{-2});
\ \ y_2(x)\sim x^{r}\exp(-\alpha\,x)(1+O(x^{-1}))\ \
(x\rightarrow\infty)
\ee

\z It is convenient to choose a    particular fundamental matrix
of the system (\ref{eqvec}) constructed with these two  solutions:

\be{choiceM}
M(x):=\pmatrix{{y_1(x)} & {y_2(x)}\cr {y'_1(x)} & {y'_2(x)}}
\sim
\pmatrix{{1+O(x^{-1})} & {(1+O(x^{-1}))x^{r}\,e^{-\alpha\,x}}\cr
         {O(x^{-2})} & {-\alpha(1+O(x^{-1}))x^{r}\,e^{-\alpha\,x}}}
\ee

\z and write the solution of (\ref{dequH}) in the form

\be{solh1}
{\bf H}(n,n+z)=M(n+z)\left(\int_{0}^{z}M^{-1}(n+t)\,{\bf T}_n(t)
  d\,t+\right. \left.
M^{-1}(n)\left({\bf H}(n-1,n_-)-{\bf E}_n\right)\right)
\ee

\z Taking $z=1$ in the relation above  we get a recurrence relation
for the $\{{\bf H}(k-1,k_-)\}_k$. With the substitution
$M^{-1}(k+1){\bf H}(k,(k+1)_-)={\bf q}_k$ it reads

\be{recubfq}
\quad\  {\bf q}_n={\bf q}_{n-1}+\int_{0}^{1}M^{-1}(n+t)\,{\bf T}_n(t)d\,t-M^{-1}(n){\bf E}_n={\bf q}_{n-1}+M^{-1}(n){\bf V}(n)
\ee
with
\be{defV(n)}
{\bf V}(n):=\int_{0}^{1}M^{-1}(n;n+t)\,{\bf T}_n(t)d\,t-{\bf E}_n
\ee

\z and where $M(n;n+t)$ is the fundamental matrix specified by $M(n,n)=I$,
whence

\be{solHn}
{\bf H}(k,(k+1)_-)=M(k+1)\left(\sum_{j=k_0+1}^{k}
M^{-1}(j){\bf V}(j)+{\bf q}_{k_0}\right)
\ee

\z For large $n$ and
$t\le 1$  we obtain in a straightforward manner,

\be{asymM}
M^{-1}(n,n+t)=
\pmatrix{1&{\frac{1}{\alpha}(1-e^{\alpha\,t})} \cr 0&{e^{\alpha\,t}}}+O(n^{-1})
\ee
\z and from (\ref{EstRem}) 
\be{approxT(n)}
\ {\bf T}_n(n+t)=(1+O(n^{-1}))\pmatrix{0\cr{-e^{-t}}}R_n(n)
\ee
Beginning with this point, most estimates will be 
different in the special case $\alpha=1$ which
corresponds to being on the  Stokes line. Combining (\ref{approxT(n)})
and Propositions {\ref{CvL}},
{\ref{RemainderEq}} we get from (\ref{defV(n)})

\be{approxdel}
\quad\ \  {\bf V}(n)=R_n(n)\left(\pmatrix{{e^{-1}}\cr -1}+\pmatrix{-1\cr 1}+O(n^{-1})\right)
=\pmatrix{{\mbox{const}}\cr {O(n^{-1})}}R_n(n)\ \ \ (\alpha=1)
\ee

\z and similarly
 
\be{approxdel2}
{\bf V}(n)={\bf O}(1)R_n(n)\ \  (\alpha\ne 1)
\ee

\z In view of (\ref{choiceM})
a direct calculation  gives for $M^{-1}(n)$

$$\pmatrix{{1+O(n^{-1})}&
{\alpha^{-1}+O(n^{-1})}
\cr{O(n^{-2})n^{-r}e^{\alpha\,n}}&{(-\alpha^{-1}+O(n^{-1}))n^{-r}e^{\alpha\,n}}}
$$
(note that the absence of exponential factors from the first row
is not an effect of an  approximation; in fact 
$\arg(\alpha)$ is arbitrary and  nothing is assumed
about the size of the exponentials).  Thus

$$
M^{-1}(n)\,{\bf V}(n)=\pmatrix{{O(1)}\cr
{O(n^{-r}e^n)}}R_n(n)=
\pmatrix{{O(n^{r-1/2}\,e^{-n})}\cr {O(n^{-1/2}\,e^{(\alpha-1)n})}}\ \
(\alpha\ne 1)
$$

\z and correspondingly

$$
M^{-1}(n)\,{\bf V}(n)=\pmatrix{\mbox{const}
\cr {\mbox{const}\,n^{-r-1}e^n}}R_n(n)=
\pmatrix{{O(n^{r-1/2}\,e^{-n})}\cr {O(n^{-3/2})}}
\ \ \ (\alpha=1)
$$

\z so that
the series $\sum_k M^{-1}(k){\bf V}(k)$ converges and 

\be{estitail1}
\sum_{k=n}^{\infty}M^{-1}(n){\bf V}(n)= 
       \pmatrix{{O(n^{r-1/2}\,e^{-n})}\cr {O(n^{-1/2}\,e^{(\alpha-1)n})}}       
\ \ \ (\alpha\ne 1)
\ee
\z and 

\be{estitail}
     \pmatrix{{O(n^{-r-1/2}\,e^{-n})}\cr {O(n^{-1/2})}}\ \ \ (\alpha=1)
\ee
\z If we make the choice 
$${\bf q}_{k_0}= -\sum_{k=k_0+1}^{\infty} M^{-1}(k){\bf V}(k)$$
in (\ref{solHn}) we get

\be{choiceHn}
{\bf H}(k,(k+1)_-)=
-M(k+1)\sum_{k+1}^{\infty}
M^{-1}(j){\bf V}(j)
\ee

\z so that in view of (\ref{estitail1}),
 (\ref{estitail}) and (\ref{choiceM})

\be{finalestrec}
{\bf H}(n,(n+1)_-)={\bf O}(n^{r-1/2}e^{-n})
\ee
for all $\alpha$ with $|\alpha|=1$
$\square$

At this point we can prove the existence part of Theorem {\ref{T1}}.
For $R\ne 0$ it is a direct consequence of the Lemma {\ref{lemmaex}},
and of the Propositions {\ref{OptSummationProp}}
and {\ref{CvL}}. 
If $R=0$ we proceed in essentially the
same way. Taking $n$ such that $x\in [n,n+1]$ we
get for the special solution provided by Lemma {\ref{lemmaex}},

\be{caseR=0}
|Y-\sum_{j=0}^{k}\frac{s_k}{x^k}|\le
|Y-Y_n|+|\sum_{j=k}^{n}\frac{s_k}{x^k}|
\ee

\z Since now we have $|s_n|<\mbox{const}\Gamma(n+\rho-1)|$ the
second sum in (\ref{caseR=0}) can be estimated 
using Proposition {\ref{uniform,factorial}} 
taking $c_j:=\mbox{const}\Gamma(n+\rho-1)$ (and $G(j):=\ln(c_j)$).
It follows immediately from (\ref{caseR=0}) that $Y$
satisfies the inequalities  (\ref{mainthm}) with $C_1=1,\ C_2=0$,
$A_1=0$ and  some $B$.
\vskip 12pt

{\em Uniqueness}.  We can assume without loss of generality (\ref{simpli})
and  $\Re(\alpha\le 0)$. The case $\Re(\alpha)=0$ is
actually trivial since it reduces to a statement about
asymptotics to all orders. We show that 
if $y$ satisfies (\ref{mainthm}) then $y=C_1 y_1+C_2 y_2$
where $y_{1,2}$ are the solutions corresponding
to $(C_1,C_2)=(1,0)$ and $(0,1)$ respectively
whose existence has already been proven.

We have to distinguish the case $A_1\ne 0$ (which implies $\alpha=1$
and $R\ne 0$, cf.(\ref{defA12})).  Taking $k_2=0$ (which means that we
only keep the leading term of the second power series) it follows from
(\ref{mainthm}) that

$$\left|F(x)-C_1\sum_{j=0}^{k}\frac{s_{1,j}}{x^j}
-C_2\,x^{r}\,e^{-x}
\right|<
$$
\vskip -0.5cm
$$|C_1|(A_1\sqrt{K_1}+B)\,\left|x^{-k-1}
\Gamma(k+r+1)\right|+
\mbox{const}\,|C_2|\left|x^{r-1}
\,e^{-x}
\right|
$$

The case $C_1=0$ is trivial 
because the only solutions that decay exponentially
at infinity are multiples of $y_2$. With $C_1\ne 0$
an elementary computation shows that one can choose a
$B'$ so that for all $x\in {\cal R}_{\theta}$ and $k\in \NN$

$$|C_1|(A_1\sqrt{k}+B)\,\left|x^{-k-1}
\Gamma(k+r+1)\right|+
\mbox{const}\,|C_2|\left|x^{r-1}
\,e^{-x}
\right|<
$$
\vskip -0.5cm
$$
|C_1|(A_1\sqrt{k}+B')\,\left|x^{-k-1}
\Gamma(k+r+1)\right|
$$

\z Consequently,
the function $y(x)-C_2\,x^{r}\,e^{-x}$ is optimally asymptotic
to the series $C_1\,S_1$ (since $R\ne 0$). 
Because  $C_1\,y_1(x)$ has the same property,
we have, in view of the corolarry to Proposition {\ref{OptSummationProp}}

$$|y(x)-C_2\,x^{r}\,e^{-x}-C_1\,y_1(x)|=o(x^{r}\,e^{-x})=o(y_2)
 {\mathrm\ as}\ x\rightarrow\infty$$

\z which means that $y(x)-C_1\,y_1(x)=C_2\,y_2(x)$.

If now $A_1=0$ we take $k_1=k_2=[x]$ in 
(\ref{mainthm}) and get in a straightforward way:

$$\left|y(x)-C_1\sum_{j=0}^{[x]}\frac{s_{1,j}}{x^j}
-C_2\,x^{r}\,e^{-x}\sum_{j=0}^{[x]}\frac{s_{2,j}}{x^j}
\right|<\mbox{const}\,x^{r-1/2}\,e^{-x}
$$

Since by Lemma {\ref{lemmaex}} the same inequality
is true with $y(x)$ replaced by $C_1 y_1(x)+C_2 y_2(x)$
it follows that $|y(x)-C_1 y_1(x)-C_2 y_2(x)|$ decays faster
than $y_2(x)$ as $x\rightarrow\infty$ which is possible
only if it vanishes identically. $\square$
\vskip 1ex
{\em{Proof}} of Lemma \ref{Puniq}.

For this purpose we only need the recurrence (\ref{recp=1}), the
asymptotic behavior of its solutions (\ref{asymptcoefLin}) and the
condition (\ref{Gevrey2}). Assume we have a formal power series
$S:=\sum_{k=0}^{\infty}s_{k}x^{-k}$ that
solves a canonical differential equation of the type considered.
The recurrence (\ref{recp=1}) provides us in a
straightforward way with a set of equations for $a_j, b_j$,
$j=0,1,2,...$:

\be{set eq 1}
-s_1 a_j+s_0 b_{j+2}=T_j \mbox{\quad for $j=0,1,2,...$}
\ee

\z where $T_j$ only depend on the series $S$ and on the coefficients
$a_i, b_{i+2}$ with $i<j$. 

The equations (\ref {set eq 1}) alone would not determine the
$a_i,b_i$ uniquely but the condition (\ref{Gevrey2}) binds $a_i$ and
$b_i$ together. The best way to see this is to divide (\ref{recp=1})
by $s_n$, use (\ref{asymptcoefLin}), (\ref{Gevrey2}) (and the
assumption $R\ne 0$) to write an asymptotic expansion to all orders in
$\frac{1}{n}$ of the resulting equation and then equate the successive
powers of $\frac{1}{n}$.

 The upshot is the system:

\be{set eq 2}
a_0 a_j-b_j=\tilde T_j \mbox{\quad for $j=2,3...$}
\ee

\z where $\tilde T_j$ depend on the series $S$ and on $a_i,b_i$  with $i<j$.
But $a_0$ and $a_1$ are determined from (\ref{asymptcoefLin}), namely
$a_0=\alpha,\ a_1=-r$, so that  (\ref{set eq 1}), (\ref{set eq 2})
determine uniquely the $a_i,b_i$. $\square$
\vskip 1ex
{\centerline{*}}
\vskip 1ex

{\em{Proof}} of Lemma \ref{lemma5}.

 For simplicity
we make a change of variables so that 

$$|\alpha|=1;\ \ x\in\RR^+$$

The recurrence relation for the formal series solution 
is obtained in the usual way, by inserting the formal series in the
differential equation, expanding everything out
and identifying the powers of $x$:

$$-\sum_{k=2}^{\infty}\frac{(k-1)s_{k-1}}{x^k}=$$
\vskip -0.5cm
$$
\sum_{j=0}^{\infty}\frac{f_{0,j}}{x^j}+
\sum_{j=0}^{\infty}\frac{f_{1,j}}{x^j}\sum_{j=1}^{\infty}\frac{s_j}{x^j}
+..+\sum_{j=0}^{\infty}\frac{f_{P,j}}{x^j}\left(\sum_{j=1}^{\infty}\frac{s_j}{x^j}
\right)^P=
$$
\be{compactf}
\sum_{k=1}^{\infty}\frac{1}{x^k}\sum_{J_0}f_{p,i}s_{k_1}...s_{k_p}
\ee
where $J_0$ consists of all the integer tuples $(i,k_1,...,k_p)$,
$k_i\ge 1$ and $i\ge 0$, with
$k_1+k_2+..+k_p+i=k$. It follows,

\be{recuformal}
s_k=\frac{1}{\alpha}(k+r-1)s_{k-1}+\sum_{J_1} f_{p,i}s_{k_1}...s_{k_p}
\ee

\z where $r=f_{1,1}+2f_{2,0}s_1$
and the index set $J_1$ in the sum excludes the tuples with $k_i>k-2$
from $J_0$. We make the substitution (suggested by solving the
linearized version of the recurrence)

$$s_k=\alpha^{-k}\Gamma_{k+r}\eta_k$$

\z where $\Gamma_{k+r}=\Gamma(k+r)$ if $\Re(k+r)>-1$ and $\Gamma_{k+r}=1$
otherwise (again, to avoid unpleasant poles). For large $k$,
\be{recufgen}
\eta_k=\, \eta_{k-1}+\sum_{J_1	} 
        {f_{p,i}\eta_{k_1}...\eta_{k_p}}\alpha^i
         \frac{\Gamma_{{k_1}+r}...\Gamma_{{k_p}+r}}{\Gamma_{k+r}}
\ee

\z The existence of a formal solution of the differential equation
follows from the existence of a solution to the recurrence relation
(\ref{recufgen}) which is obvious. The first objective is then to show
that the sequence $\{\eta_k\}_k$ converges; the only delicate
step is to show that  the sequence $\{|\eta_k|\}_k$ is bounded.
To this end we compare the recurrence (\ref{recufgen})
with a suitable {\em linear} recurrence (\ref{reculin}) below.

Taking out of the sum the terms with $k_i>k-n$ for a conveniently
large $n$ we get, for some constants $A_1,..,A_n$,

\be{recufnp}
\eta_k=\eta_{k-1}+ \frac{A_1}{k^2}\eta_{k-2}+ 
..+\frac{A_n}{k^n}\eta_{k-n}+\sum_{J_2} 
        {f_{p,i}\eta_{k_1}...\eta_{k_p}}\alpha^i
         \frac{\Gamma_{{k_1}+r}...\Gamma_{{k_p}+r}}{\Gamma_{k+r}}
\ee

\z where the restriction $J_2$ in the sum now reads $k_1+k_2+..+k_p+i=k; k_i<
k-n$. Taking the initial condition ${\overline{\eta_1}}=|\eta_1|$
in the recurrence

\be{recufnpmod}
\ \ \ \ \ \overline\eta_k=\overline\eta_{k-1}+ \frac{|A_1|}{k^2}\overline\eta_{k-2}+ 
..+\frac{|A_n|}{k^n}\overline\eta_{k-n}+\sum_{J_2} 
        {f_{p,i}\overline\eta_{k_1}...\overline\eta_{k_p}}
         \frac{\Gamma_{{k_1}+\rho}...\Gamma_{{k_p}+\rho}}{\Gamma_{k+\rho}}
\ ( \rho=\Re(r))
\ee

\z we clearly get
$\overline\eta_k\ge|\eta_k|$.  Now the solution of
(\ref{recufnpmod}) with positive initial conditions is increasing in
$k$ so that for $A>\max\{A_1,..,A_n\}$

\be{recn3}
\overline\eta_k\le\, \left(1+\frac{A}{k^2}\right)
        \overline\eta_{k-1}+ \sum_{J_2} 
        {f_{p,i}\overline\eta_{k_1}...\overline\eta_{k_p}}
         \frac{\Gamma_{{k_1}+\rho}...\Gamma_{{k_p}+\rho}}{\Gamma_{k+\rho}}
\ee
Taking
$\overline\eta_k=\sigma_k\prod_1^k(1+A/j^2)/\prod_1^k(1-1/j^2)$ we get for some
constant $C$
\be{recuq}
\sigma_k\le\, \left(1-\frac{1}{k^2}\right)
        \sigma_{k-1}+ C\sum_{J_2} 
        {f_{p,i}\sigma_{k_1}...\sigma_{k_p}}
        \frac{\Gamma_{{k_1}+\rho}...\Gamma_{{k_p}+\rho}}{\Gamma_{k+\rho}} 
\ee

\z Noting that  card($J_2$)$\le P\,k^p$ 
we get

\be{ren1}
\sigma_k\le\left(1-\frac{1}{k^2}\right)\sigma_{k-1}+
\frac{C\,P}{k^2}\max_{J_2}\left\{k^{p+2}\,f_{p,i}\sigma_{k_1}...\sigma_{k_p}         
         \frac{\Gamma_{{k_1}+\rho}...\Gamma_{{k_p}+\rho}}{\Gamma_{k+\rho}}
\right\}
\ee

It is easy to see by a straightforward calculation that one can choose
$k$ and $n$ large enough  so  that 

\be{cond k}
CP\,k^{p+2}f_{p,i}\frac{\Gamma_{{k_1}+\rho}...\Gamma_{{k_p}+\rho}}{\Gamma_{k+\rho}}
<\left(\frac{k_1!..k_p!i!}{(k_1+..+k_p+i)!}\right)^{\frac{1}{2}}<
\left(\frac{k_1!..k_p!}{(k_1+..+k_p)!}\right)^{\frac{1}{2}}
\ee

\z when  $k_i<k-n$ (for the second inequality above note that the
middle term is a decreasing function of its arguments). Using
(\ref{ineqf_i}) we get

\be{ineq2}
\sigma_k\le\left(1-\frac{1}{k^2}\right)\sigma_{k-1}
    +\frac{1}{k^2}\max_{J_3}\{\sigma_{k_1}...\sigma_{k_p}
           \left(\frac{k_1!..k_p!}{(k_1+..+k_p)!}\right)^{\frac{1}{2}}\}
\ee

\z $J_3$ consists of all the tuples $k_1,..,k_p$ with $k_1+..+k_p\le k$.
By construction the  empty tuple contributes with  
const$\,\kappa^k\le$ const. It
follows easily that

\be{reccompn}
\sigma_k\le\max_{J_3}\{\sigma_{k-1};\sigma_{k_1}...\sigma_{k_p}
       \left(\frac{k_1!..k_p!}{(k_1+..+k_p)!}\right)^{\frac{1}{2}}\}
\ee

\z from which, if $\sigma_1$ is large enough which we are allowed to 
assume since the solution is increasing in $\sigma_1$, it follows
that $\sigma_k$ is bounded by the solution of the recurrence:
\be{recuit}
\sigma_k=\max_{J_3}\{\sigma_{k_1}...\sigma_{k_p}
       \left(\frac{k_1!..k_p!}{(k_1+..+k_p)!}\right)^{1/2}\}
\ee

\z For each $k$, (\ref{recuit}) determines
a tuple $k^*_1,..,k^*_p,\ k^*_1+..+k^*_p\le k$ for which
the maximum is realized. If we  take $\sigma_k=\exp(\phi_k)$
we get a linear recurrence:

\be{reculin}
\phi_k=\phi_{k^*_1}+..+\phi_{k^*_p}+
           \frac{1}{2}\ln\left(\frac{k^*_1!..k^*_p!}{(k^*_1+..+k^*_p)!}\right)
\ee
 Consider the associated homogeneous recurrence
\be{recuiteg}
F_k=F_{k^*_1}+...+F_{k^*_p};\ \ F_1=1
\ee

\z Inductively we see that:

\be{estimhm}
F_k\le k
\ee
We claim that 
\be{estimin}
\phi_k\le \phi_1\,F_k-\frac{1}{2}\ln(F_k!)
\ee

\z Indeed, the inequality above is true for $k=1$ and
assuming it holds for all $k'<k$ we have (cf. also the comment
following (\ref{cond k}))

$$\phi_k\le\phi_1(F_{k^*_1}+..+F_{k^*_p})
-\frac{1}{2}\ln(F_{k^*_1}!..F_{k^*_p}!)
+\frac{1}{2}\ln(\frac{F_{k^*_1}!..F_{k^*_p}!}{(F_{k^*_1}+..+F_{k^*_p})!})$$
in view of (\ref{estimhm}) which gives (\ref{estimin}).

But note that the function
$g(x)=A\,x-B\,\ln(x!):(1,\infty)\mapsto\RR$ is bounded from above so
that the solution of (\ref{recufgen}) is bounded. 

Now, taking $n$ large
enough one can make the nonlinear term in (\ref{recufnp}) smaller than
$\mbox{const}\,k^{-m}$ (cf. also (\ref{cond k})) with $m$ as large needed.
The rest of the proof is obvious. $\square$

We return to the proof of Theorem \ref{T2}.

{\em{Proof.}} We start
by showing the existence of a special solution corresponding
to $C=0$.  We look for solutions of the equation (\ref{1stnonln}) in 
the form
\be{ansatz1}
y(x)=Y_n(x)+V(n;x);\ \ Y_n:=\sum_{k=1}^n\frac{s_k}{x^k};\ \ (x\in [n,n+1])
\ee

\z where $n=[x]$. The differential equation for  $V$ reads:
\be{deV}
V'=\frac{d}{dy}F(Y_n,x)\,V+\sum_{p=2}^{P}\frac{d^p}{dy^p}F(Y_n,x)\,V^p
   +T_n
\ee

\z where 

$$T_n=F(Y_n,x)-Y'_n(x)$$
to which we add the condition of  continuity of $y$:

\be{conV}
V(n,n_+)=V(n-1,n_-)-\frac{s_n}{n^n}
\ee

\z \bl{lemmap}
There exists a solution of the problem (\ref{deV})--(\ref{conV})
with the property

$$|V(n,x)|\le \mbox{const}\,\frac{\Gamma_{n+\rho}}{x^n}\ \ (x\in [n,n+1])$$
\end{Le}

\z Again we need first some estimates on the power series solution.

\z \bp{estimT}
For $x\in [n,n+1]$

$$T_n=(R+O(\frac{1}{n}))\alpha^{-n-1}\frac{\Gamma_{n+\rho+1}}{x^{n+1}}$$

\end{Pro}
{\em Proof

}

\z For $p\in\{0,..,P\}$ we write 

\be{condFi}
F_p(x)=\sum_{k=0}^n\frac{f_{p,k}}{x^k}+\frac{R_{p,n}(x)}{x^{n+1}}
\mbox{\ \ with\ \ } {R_{p,n}}(x)<C_F\kappa^n n!  \ \ (x\in[n,n+1])         
\ee

\z Then $T_n$ is given by
\be{dif}
\sum_{p=0}^P Y_n^p\left(\sum_{k=0}^n\frac{F_{p,k}}{x^k}
+\frac{R_{p,n}}{x^{n+1}}\right)-Y_n'=
\ee
\vskip -0.3cm
$$\sum_{Q}\frac{1}{x^{Q}}\sum_{J_Q}
f_{p,i}s_{k_1}...s_{k_p}+\sum_{k=1}^{n}\frac{k\,s_k}{x^{k+1}}+
\sum_{i=1}^P\frac{R_{p,n}}{x^{n+1}}Y_n^p
$$

\z where the index set $J_Q$ contains tuples $(i,k_1,..,k_p)\in\{0,..,n\}$ with
 $i+k_1+..+k_p=Q$.  The recurrence relation for $s_k$ is such
that all the terms with $Q\le n$ have to compensate each other whereas
for $Q=n+1$ they add up to $\alpha {s_{n+1}}{x^{-n-1}}$. We thus get for
$T_n$

\be{T_n}
\left|T_n-\alpha\frac{s_{n+1}}{x^{n+1}}\right|\le
\sum_{Q\ge n+2}\frac{1}{x^{Q}}\sum_{J_Q}
|f_{p,i}||s_{k_1}|...|s_{k_p}|+
\sum_{i=1}^P\frac{|R_{p,n}|}{x^{n+1}}|Y_n^p|
\ee

The last term on the rhs of (\ref{T_n}) can be estimated
by $\mbox{const}\,\Gamma_{n+\rho+1}\,x^{-n-2}$. For the first term we write

\be{estimsum}
\left( \sum_{n+2\le Q\le n+s}+\sum_{n+s\le Q\le (P+1)n}\right)\frac{1}{x^{Q}}\sum_{J_Q}
|f_{p,i}||s_{k_1}|...|s_{k_p}|
\ee

\z For all $n$, $|s_n|$ and $|f_{p,n}|$ are bounded by $\mbox{const}\,
\Gamma_{n+\rho}$. Note that the $\Gamma$ function is log-convex so
that the maximum of a product of the form $\Gamma(x_1)..\Gamma(x_m)$
over a convex domain, is reached on the boundary of the domain. Now,
if $n$ is large enough compared to $s+p$,

$$
\sum_{\stackrel{{ {\scriptstyle k_1+..+k_p=n+s}}}{k_i\le
n}}\Gamma_{k_1+\rho}..\Gamma_{k_p+\rho} \le
\mbox{const}\, {P\sum_{j=0}^{p+s+1}n^{-j}
\sum_{\stackrel{{ {\scriptstyle 
k_1+..+k_{p-1}=s+j}}}{ k_i\le
n}}\Gamma_{k_1+\rho}..\Gamma_{k_{p-1}+\rho}}  +
$$
\be{elemineq1}
\mbox{const}\,{{n+p+s}\choose{p}}n^{-p-s-2}<
\mbox{const}(P,s)\Gamma_{n+\rho}
\ee

\z It follows that  the leftmost sum in (\ref{estimsum})
is less than $\mbox{const}\,{\Gamma_{n+\rho}}{x^{-n-2}}$ 
for $x\in[n,n+1]$. In the second sum we write
 $$Q=q\,n+m$$

\z with $q\le P, m\le n$. The log-convexity argument
together with the restriction on the indices $i, k_j$ give

$$
|f_{p,i}||s_{k_1}|...|s_{k_p}|<\mbox{const}\,\kappa^i
\,i!\Gamma_{k_1+\rho}..\Gamma_{k_p+\rho}\le
$$
\vskip -0.3cm

\be{estsum_3}
\mbox{const}\,(\Gamma_{n+\rho})^q
\Gamma_{m+\rho}
\ee

\z and we crudely bound the second sum in (\ref{estimsum}) by the number
of terms times the maximal term. An elementary application of
Stirling's formula gives

\be{estis2}
\mbox{const} \,n^{P+2} \max_{n+s\le Q\le nP}\left\{
          \frac{\Gamma_{n+\rho}^q\Gamma_{m+\rho}}{n^{qn+m}}
\right\}
\le
\ee
\vskip -0.3cm
$$\mbox{const}\,n^{P+2}\exp\left(\max_{n+s\le Q\le
nP}\left\{(m+\rho-\frac{1}{2})\ln(m)
+(q\rho-\frac{q}{2}-m)\ln(n)-m-qn\right\}\right)$$

The expression to be maximized is decreasing in $q$ and is convex in
$m$. It is not difficult to see that for $n$ large enough the maximum
is reached at $q=1,m=s$.  For large $n$ (\ref{estis2}) is bounded by
$\mbox{const}\,n^{r-s-1/2+P+2}\,e^{-n}$. Choosing $s>P+4$ we make the second
sum in (\ref{estimsum}) less than ${\Gamma_{n+\rho}}{n^{-n-2}}$. The
proof is completed by combining these inequalities with (\ref{T_n})
$\square$

We take now $n$ large enough and first find a suitable solution of the
linearized version of the equation (\ref{deV}) on each interval
$[n,n+1]$. Consider the differential equations:

\be{delH}
H'=F_y(Y_n,x)H+T_n
\ee
with the initial conditions at $x=n$ 
chosen so that $y(x)$ as defined in (\ref{ansatz1})
is continuous namely,
\be{condcap}
H(n,n_+)=H(n-1,n_-)-\frac{s_n}{n^n}
\ee
under the same hypothesis on $F$ as before.
\bp{1stnlex}
There is a solution to the problem (\ref{delH}), (\ref{condcap})
such that 
\be{estim Hglob}
|H(n,x)|\le \mbox{const}\,\frac{\Gamma_{n+\rho}}{x^n}\ \ (x\in [n,n+1])
\ee
\end{Pro}

\vskip 1ex

\z The proof is the one-dimensional projection of the proof
 of lemma {\ref{lemmaex}}\ $\square$.

It also follows that for $x\in [n,n+1]$

\be{estimerror}
H(n,x)\le \mbox{const}\,x^{\rho-1/2}\,e^{-x}
\ee

\vskip 2ex

We now turn to the proof of Lemma \ref{lemmap}. We shall look for
solutions of the equation (\ref{deV}) in the form $V=H+\delta$ where
$H$ is the function provided by the Proposition \ref{1stnlex}. Then
$\delta$ is a smooth function and satisfies the differential equation

\be{equdelta}
\delta'=\frac{d}{dy}F(Y_n,x)\,\delta+
\sum_{p=2}^{P}\frac{d^p}{dy^p}F(Y_n,x)\,(H+\delta)^p
\ee

We are looking for a solution of (\ref{equdelta}) of the order
$O(H^2)$; such a solution is constructed as a fixed point of a
contractive mapping.

Take $n$ large enough, $\beta>2\rho-1$ and let ${\cal{S}}_{n,\beta}$ be the
space of continuous functions on $(n,\infty)$ with the norm
\be{defnorm}
\|f\|:=\sup_{t>n}\left|t^{-\beta}e^{2\,t}\,f(t)\right|
\ee

and the operator 
\be{def L}
({\cal L}\, f)(x):=\int_x^{\infty}
\frac{d}{dy}F(Y_n,t)\,f(t)+
\sum_{p=2}^{P}\frac{d^p}{dy^p}F(Y_n,t)\,(H(t)+f(t))^p\,d\,t		
\ee

In view of (\ref{estimerror}) for large $n$,
$H^2\in{\cal{S}}_{n,\beta}$.  Actually we can make the norm of $H^2$ very
small by choosing $n$ large enough.  ${\cal L}$ maps
${\cal{S}}_{n,\beta}$ into itself and for a given  $\epsilon$ there is
$n_{\epsilon}$ large enough so that ${\cal L}$ maps ${\cal B}_{\epsilon}$,
the ball of radius $\epsilon$, into itself. Actually, ${\cal L}$ is a
contraction for large $n$. Indeed

\be{estimop}
{\cal L}\, f=\int_x^{\infty}\alpha \,f(t)d\,t+O(f/x,f^2,H\,f,H^2)
\ee

\z where $O$ is with respect to the norm (\ref{defnorm}).  But if
$f,g\in{\cal B}_{\epsilon}$ and 
$\epsilon'>0$ it is straightforward to see that 
$\|{\cal L}f-{\cal L}g\|<
\frac{1+\epsilon'}{2}\|f-g\|$ if $n$ is large in which case 
the equation (\ref{equdelta}) has a solution in ${\cal{S}}_{n,\beta}$.
It follows that the equation (\ref{deV}) has a solution
with the required properties $\square$.

Lemma {\ref{lemmap}} together with Proposition
{\ref{OptSummationProp}} prove the existence part of Theorem
{\ref{T2}} for $C=0$. Given this, for $C\ne 0$ the existence result
follows easily, since for any $C$ there exists a unique solution of
the differential equation (\ref{1stnonln}) of the form
$y(x)=V(x)+\epsilon(x)$ with $V(x)$ the solution given in Lemma
{\ref{lemmap}} and $\epsilon(x)\sim C e^{-\alpha x}(x^{r+1}+O(x^{r}))\
(|x|\rightarrow\infty)$. This is a standard result (see \cite{Wasow},
\cite{Coddington}) ; in our case $\epsilon(x)$ could be also directly
obtained using the contractive mapping arguments above.

The uniqueness proof is very similar to the one in the linear case
since no two solutions can differ by less than $\mbox{const}\,e^{-\alpha
x}x^{\rho+1}$ as $x\rightarrow\infty$ without being equal to each
other. Finally, the information that a solution $y(x)$ is decaying at
infinity is enough to guarantee that $y$ is asymptotic to the series
(\ref{defseri}).  The difference $y(x)-V(x)$ has the asymptotic
behavior $C e^{-\alpha x}(x^{r+1}+O(x^{r}))$ for some $C$. Let $y_C$
be the special solution satisfying (\ref{eqT2}). Since $y(x)-V(x)\sim
C e^{-\alpha x}(x^{r+1}+O(x^{r}))$ as well, it follows that
$y(x)-y_C(x)=o(e^{-\alpha x}x^{r+1})$ which means $y(x)=y_C(x)$
(\cite{Coddington}; a simple direct proof uses contraction mapping
(\ref{def L})). The proof of Theorem \ref{T2} is complete. $\square$

\vskip 2ex
{\centerline{*}}
\vskip 2ex

{\em Discussion of the examples}.

\z  i) Proof of Eq. (\ref{eqEi}). We are looking
for a constant $C$, dependent on the ray in the complex plane,
such that

\be{aimEi}
E_n(x):={\cal E}i(x)-e^x\sum_{k=1}^n \frac{(k-1)!}{x^k}=C+o(1)
\ (|x|\rightarrow\infty, \ \arg(x)=\theta)
\ee

\z By Theorem 1.1 and Proposition \ref{OptSummationProp}, 
there exists such a constant. We might
as well compute it for $|x|\rightarrow\infty$ along a subsequence,
say $x_n=ne^{i\theta}$.

$\bullet$ Take first $\theta>0$ ($\theta<0$ is  similar). It is easy 
to show that

$$E_n(x)=\pi i+
n!n^{-n}\int_C dt\frac{e^{nt}}{t^{n+1}}
$$

\z where $C$ is a contour joining $-\infty$ to $e^{i\theta}$
above the real axis; we choose it to be the stationary phase (or, which
is the same, the steepest ascent) contour for the integrand:
$\Im(nt-(n+1)\ln(t))=\mbox{const}$.


The leading behavior of the integral for large $n$ is due to the contribution
of a region, near the right end point of $C$, where $z:=t-e^{i\theta}$
is of the order $O(1/n)$.

$$\int_C dt\frac{e^{nt}}{t^{n+1}}=
(1+o(1))e^{ne^{i\theta}-(n+1)i\theta}\int_{-\infty}^0dz e^{n(1-e^{-i\theta})z}$$

\z so that, for $\theta>0$,

\be{estiEi,away}
E_n-\pi i=(1+o(1))\sqrt{\frac{2\pi}{n}}\frac{1}{1-e^{-i\theta}}
e^{n(e^{i\theta}-1)-(n+1)i\theta}=o(1)
\ee

$\bullet$ For $\theta=0$  we have:

$$E_n(n)=n!n^{-n}{\mathrm
PV}\int_{-\infty}^1 dt\frac{e^{nt}}{t^{n+1}}=
-n!n^{-n}e^n\Im\int_0^{\infty}\frac{e^{-iny}}{(1-iy)^{n+1}}dt $$ 

\z where  the last expression is obtained after pushing
the contour parallel to the imaginary axis and taking
$t=1+iy$. Now, choosing $\beta\in(-\frac{1}{2},-\frac{1}{3})$
we get

\ba{estlt}
\int_0^{\infty}\frac{e^{-iny}}{(1-iy)^{n+1}}dt&=&
\int_0^{n^{\beta}}\frac{e^{-iny}}{(1-iy)^{n+1}}dt+
O(\exp(-\frac{1}{2}n^{2\beta+1}))=\cr
&&\int_0^{\infty}e^{-n\frac{y^2}{2}}\left(1+iy-\frac{i}{3}ny^3\right)dy
+o(n^{-1})=\sqrt{\frac{\pi}{2n}}+\frac{i}{3n}+o(n^{-1})
\end{eqnarray}

\z so that, after taking the imaginary part,
\be{estimfinEi} 
E_n(n)=
-\frac{1}{3}\sqrt{\frac{2\pi}{n}}(1+o(1))
\ee$\square$

For the exponential integral Ei($x$), it is easy to evaluate the
optimal weight from (\ref{estimfinEi}) and its differential equation.
The arguments below can be easily made rigorous but we will not insist on
that, since we are only aiming at an illustration.
Let 

\be{eig(n)start}
R(x)=e^{-x}\mbox{Ei}(x)-\sum_{k=0}^{n-2}\frac{k!}{x^{k+1}}
\ee

\z Then, $R(x)$ satisfies the differential equation
$$R'+R=\frac{(n-1)!}{x^n}$$

\z Taking $R(x)=\frac{(n-1)!}{x^n}g(x)$ we get for $g$ the equation
$g'+(1-nx^{-1})g=1$.
The relevant region is $x=n+O(\sqrt{n})$ 
so it is convenient to change  variables further to  $x=n+s\sqrt{n}$ 
and $g(x)=\sqrt{n}\phi(s)$. We get

$$\phi'+\frac{s}{1+sn^{-1/2}}\phi=1$$

\z with the initial condition, coming from (\ref{estimfinEi})
 $\phi(0)=O(n^{-1/2})$. The solution is given by

$$\phi(s)=e^{-\frac{s^2}{2}}\int_0^s e^{\frac{t^2}{2}}dt+O(n^{-1/2})$$

\z The maximum value of $|\phi(s)|$ is therefore equal to $a_*$.

\vskip 0.3cm 
ii) Derivation of (\ref{AsyAirygen}).  To make the calculation of the
asymptotic series easier we take further $g(s)=s^{-1/6}h(s)$ (to make
$r_1=0$) and finally $s=3t/4$ (to have $\lambda_1=0;
\lambda_2=-1$). The resulting equation  is

$$h''+h'+\frac{5}{36t^2}h=0$$

\z having the following general formal solution:

$$C_1\sum_{j=0}^{\infty}\frac{h_j}{t^j}+C_2e^{-t}
\sum_{j=0}^{\infty}\frac{(-)^j h_j}{t^j}$$

\z The recurrence for $h_j$:
 $$(j+1)h_{j+1}=(j+\frac{1}{6})(j+\frac{5}{6})h_j$$

\z  with the choice  $h_0=1$, gives

$$h_n=\frac{\Gamma(n+\frac{5}{6})
\Gamma(n+\frac{1}{6})}{2\pi\Gamma(n+1)}
=\frac{1}{2\pi}\Gamma(n) \left(1+O(\frac{1}{n})\right)
\ \ \mbox{for large $n$}
$$

\z so that with this choice we get $R_1=R_2=\frac{1}{2\pi}$ and
$f_k=\left(\frac{4}{3}\right)^{\frac{2k}{3}}h_k$. The contents  
of (\ref{AsyAirygen}) is given by Theorem\,\ref{T1}, after making 
the appropriate  substitutions and change of variable back into  
eq. (\ref{mainthm}).

\end{section}

\begin{section}{Convergence of  asymptotic series 
 and a decomposition property}

In this section we give an example of a case of non-generic
convergence of the asymptotic series and discuss its relevance for
asymptotics beyond all orders.  The Proposition below could be easily
generalized to a larger class of differential equations but we now
look, for the sake of simplicity, at equations of the form

\be{dedeco}
\psi'+\psi=R(x)
\ee

\z where $R(x)$ is a rational function. We consider (\ref{dedeco}) for
large $x$ on a ray $\arg(x)=\theta\in(-\frac{\pi}{2},\frac{\pi}{2})$.
If $R(x)$ is a polynomial the equation has an explicit solution of the
form {\em Polynomial(x)} $+Ce^{-x}$ so the interesting case is

\be{ratdec}
R(x)=\sum_{k,j=1}^{m}R_{k,j}(x-r_j)^{-k}
\ee

\bp{deco}
\z i) Given $R(x)$ of the form (\ref{ratdec}) 
there exists a unique constant $K$ such that the solution of the
equation

\be{dedecohol}
\psi'+\psi=R(x)-\frac{K}{x}
\ee

\z is holomorphic in a neighborhood of infinity. Consequently,
any solution of the equation (\ref{dedeco})
has a decomposition:

\be{lemdeco}
\psi(x)=K\,{e^{-x}\cal E}i(x)+H\left(\frac{1}{x}\right)+C\,e^{-x}
\ee

\z where the constants $K,C$ as well as the analytic 
function  $H(\cdot)$ are uniquely determined.

\z ii) Let $\theta=0$. 
A solution of (\ref{dedeco}) has (in the sense of Theorem {\ref{T2}})
the asymptotic representation

\be{lemdecoasy}
\psi(x)\simeq \sum_{j=0}^{\infty}\frac{s_{k}}{x^k}+C\,e^{-x}
\ee

\z Moreover the constant $C$ is the same as in (\ref{dedeco})
and the power series is (factorially) divergent unless $K=0$.

\end{Pro}
{\em{Comments}}.  On the one hand the proposition above indicates in
what sense divergence of the asymptotic series is generic. On the
other hand the decomposition (\ref{lemdeco}) suggests another point of
view on the problem of the terms beyond all orders for the equation
(\ref{dedeco}). Since in any reasonable definition of asymptotic
representations, a function which is analytic at infinity should be
represented by its own (convergent) asymptotic series, once the terms
beyond all orders for a particular function (the exponential intergal)
are defined, they can be determined unambiguously for the solutions of
(\ref{dedeco}) with any rational inhomogeneity $R$. Part ii) shows
that the results obtained in this way are consistent with those
obtained through asymptotic estimates.

{\em Proof

} We have

\be{lapdec}
R(x)=\int_0^{\infty}e^{-t\,x}
           \sum_{k,j=1}^{m}R_{k,j}e^{r_j\,t}\frac{t^{k-1}}{(k-1)!}d\,t
\ee

\z Let 

\be{val1}
K=\sum_{k,j=1}^{m}R_{k,j}e^{r_j}\frac{1}{(k-1)!}
\ee

\z Then the function:

\be{complp}
\tilde R(t):=\frac{1}{1-t}\left(
\sum_{k,j=1}^{m}R_{k,j}e^{r_j\,t}\frac{t^{k-1}}{(k-1)!}-K
\right)
\ee

\z is entire so that the function

\be{analap}
z\int_{0}^{\infty}e^{-t}\tilde R(z\,t)d\,t
\ee

\z is analytic in $z$  for $|z|<\frac{1}{\rho}$. Indeed, the
integrand is analytic in $z$ and,
 $|\tilde R(z\,t)|<z^{m}t^m\,e^{|z|\,\rho\,t}$ so that
the integral is uniformly convergent for $|z\,\rho|<1$.
Furthermore, the function:

\be{sollap}
\frac{1}{x}\int_{0}^{\infty}e^{-t}\tilde R\left(\frac{t}{x}\right)d\,t+K{\cal E}i(x)
\ee

\z is a solution of the equation

\be{soleq}
\psi'+\psi=R(x)-\frac{K}{x}
\ee

\z as it can be easily checked. 
The rest of the proof of {\em i)} is immediate. 

Let 

$$H(x)=\sum_{k=0}^{\infty}\frac{h_k}{x^k}$$

\z Since clearly $\frac{\scriptstyle h_n}{\scriptstyle n!}
\rightarrow 0$ as $n\rightarrow\infty$
the proof of {\em ii)} is an easy  application of Proposition 
{\ref{OptSummationProp}} and of formula (\ref{eqEi}).$\square$

\end{section}

\begin{section}{Acknowledgments} One of the authors (O.C.) would like
to thank Professors Michael Berry, Pavel Bleher and Antti Kupiainen
for very interesting discussions. Special thanks are due to 
Prof. Joel Lebowitz for his
caring support and encouragements throughout this work.

\end{section}

\end{document}